\documentclass{amsart}

\newtheorem{theorem}{Theorem}
\newtheorem{lemma}[theorem]{Lemma}

\begin{document}
\title{Boundary layers and the vanishing viscosity limit for incompressible 2D flow}
\author{Milton C. Lopes Filho}
\date{\today}
\newcommand{\vare}{\varepsilon}
\newcommand{\real}{\mathbb{R}}
\newcommand{\disk}{\mathcal{D}}
\newcommand{\loc}{\mbox{loc}}

\begin{abstract}
This manuscript is a survey on results related to boundary layers
and the vanishing viscosity limit for incompressible flow. It is the
lecture notes for a 10 hour minicourse given at the Morningside
Center, Academia Sinica, Beijing, PRC from 11/28 to 12/07, 2007. The
main topics covered are: a derivation of Prandtl's boundary layer
equation; an outline of the rigorous theory of Prandtl's equation,
without proofs; Kato's criterion for the vanishing viscosity limit;
the vanishing viscosity limit with Navier friction condition;
rigorous boundary layer theory for the Navier friction condition and
boundary layers for flows in a rotating cylinder.
\end{abstract}

\maketitle

\section{Introduction}

    In 1904, the issue of heavier-than-air, self propelled flight by
human-made machines was at the very edge of both science and
technology. The first such flight, by the Wright brothers, occurred
at December 14th, 1903. A Brazilian author is honor bound to remark
that a more satisfying, and better publicized ``first flight'' was
achieved by Santos Dumont, a Brazilian living in France, in
September, 1906. The flight of a fixed-wing airplane could, at least
in principle, be described by near-steady, zero-viscosity,
irrotational theory of airfoils, which was already available at the
beginning of the twentieth century.

    Classical airfoil theory explained satisfactorily the
balance of forces in a wing in steady flight. In short, the force
that air exerts on the wing is divided into two standard components:
the lift (vertical force) and the drag (horizontal force), where
horizontal means the direction of steady motion. In steady flight,
these forces are balanced by the weight and the propulsion force.
The theory predicts that the lift and the drag are proportional to
the {\it circulation} of air velocity around the airfoil, and it was
in agreement with experiments, see \cite{something} for details.
Trouble occurs when one wants to change the lift, as one should do
when attempting to take off or land in a fixed wing aircraft. A
theorem, due to Lord Kelvin, states that the circulation around a
material curve, such as the boundary of an airfoil, is a constant of
motion in ideal (i.e. non-viscous) flow - or maybe, nearly constant
in slightly viscous flow. So, changing propeller speed and moving
control surfaces does not change the circulation. Since airplanes
start out at rest, with zero circulation around the wings, no
airplane could, on theoretical ground, develop a lift, and therefore
fly. Something was clearly wrong with the theory.

    The correction was due to the young theoretical mechanician Ludwig
Prandtl (1875 - 1953), who published a short paper in the
Proceedings of the Third ICM (Heidelberg 1904) whose German title
roughly translates as ``fluid flow in very little friction''. In
this article, Prandtl established a perfectly satisfactory, and
revolutionary, explanation of the following observation:

\begin{itemize}
\item[(O)]{\em The interaction of incompressible flow with a material
boundary is completely different if the flow has very small
viscosity or none at all.}
\end{itemize}

This observation, the associated explanation, called {\it boundary
layer theory} and some of what mathematicians made of this subject
in the following century and a bit, make up the subject of these
lectures. For a thorough account of the development and
understanding of the physics of boundary layers, we would like to
refer the reader to the classical book \cite{Schlichting}. It can be
argued that this short paper by Prandtl marks the birth of modern
applied mathematics.

The theory of boundary layers is a cornerstone of modern fluid
mechanics, but, as in much of this field, it lacks a {\it rational}
framework, i.e. a rigorously established connection with first
principles. Although substantial mathematical work has been done in
this direction, some basic questions remains unanswered. The purpose
of these lectures is to probe the boundaries in the mathematical
understanding of the interaction between nearly ideal flow and solid
objects, perhaps to bring what is not known about this question more
sharply into focus. The choice of material covered is strongly
slanted towards recent work by the author and his collaborators, and
it includes detailed consideration of Kato's criterion for the
vanishing viscosity limit in a bounded domain, a long discussion on
the vanishing viscosity limit for incompressible flows with Navier
boundary condition and the detailed behavior of circularly symmetric
flow inside a rotating cylinder. The choice of working with two
dimensional flow is both a reasonable pedagogical choice and a
comfort zone for the author - in the issue of boundary layers, the
sharp distinction in behavior between 2D and 3D flows is not yet
apparent, and much of the work we will discuss here generalizes
readily to 3D. Finally, we mention that these notes are written
thinking of a reasonably mature audience - we assume, not only
familiarity with standard PDE theory, but some familiarity with the
basics of mathematical fluid dynamics as well.

The remainder of these notes is divided in seven sections as
follows: Section 2 contains a derivation of Prandtl's equation;
Section 3 contains a broad overview of rigorous results on Prandtl's
equation, including some of O. Oleinik's work, and more recent
progress; Section 4 introduces and proves Kato's criterion and some
related results; Section 5 is concerned with vanishing viscosity
under Navier friction conditions and a proof of $L^p$ vorticity
estimates in this case; Section 6 contains an exposition on a
rigorous method to treat boundary layer expansions based on ideas of
geometric optics, applied to the Navier boundary condition; Section
7 explores a nearly explicit example of the behavior of the boundary
layer for the no-slip condition; Section 8 contains some conclusions
and open problems.

\section{Prandtl's theory}

In this Section, we present an asymptotic derivation of Prandtl's
boundary layer equation. Our point of departure is the Navier-Stokes
equations, which are an expression of Newton's second law applied to
the motion of a fluid, subject to an incompressibility constraint.
We write
\begin{equation} \label{nseq}
\left\{  \begin{array}{l}
\partial_t u + u \cdot \nabla u = - \nabla p + \mu \Delta u \\
\mbox{div }u = 0, \end{array} \right.
\end{equation}
where $u$ is the fluid velocity, $p$ is the scalar pressure and
$\mu$ is the kinematic viscosity of the fluid. We have assumed that
mass units have been chosen so that fluid density is one.

For this derivation, and for most of the discussion in these
lectures, we will assume that we are discussing a two-dimensional
fluid occupying a half plane $\mathbb{H} \equiv \{x_2 > 0\}$. Two
dimensional fluids really occur either in a computer simulation, or
as three dimensional fluids which are translation-invariant along
some direction (often an unstable state of affairs) or as an
approximate model for thin fluid layers. Much of the discussion
extends very naturally to a fairly general fluid domain in three
dimensions, but we will stay with the simplest possible situation
for pedagogical reasons.

As we all know, problem \eqref{nseq} requires a boundary condition
at $\{x_2 = 0\}$, and the condition usually deemed appropriate is
the {\it no slip} condition $u(x_1,0,t) = 0$. This condition
expresses an assumption that viscous fluids adhere to material
objects, something that is neither physically nor mathematically
obvious and was subject for heated debate until the mid nineteenth
century, when it became clear that it gave good agreement with
experiments.

Ideal, or inviscid flow is represented by solutions of Euler`s
equations, which is system \eqref{nseq} with $\mu = 0$. For ideal
flow, the correct boundary condition is the non-penetration
condition $u_2(x_1,0,t) = 0$. The Navier-Stokes system is a singular
perturbation of the Euler system, because the small constant $\mu$
appears in front of the highest order term of \eqref{nseq}. One
consequence of this singular perturbation is the disparity in
boundary conditions between Euler and Navier-Stokes flows - namely
that $u_1$ at the boundary goes from being identically zero for any
positive viscosity to some (in principle) nonzero function when $\mu
= 0$. This disparity is the root cause of the boundary layer
trouble.

Our objective here is to derive Prandtl`s boundary layer equations.
This is a way to quantify nearly ideal fluid behavior near the
boundary by means of an appropriate set of limit equations. Let us
begin with a simplifying assumption: we assume that the disparity
between ideal and viscous flow is concentrated on a thin layer near
$x_2 = 0$.

Our next step is to non-dimensionalize equation \eqref{nseq}, using
the time scale $T$, the length scale $L$ for horizontal lengths, and
a reference vertical length scale $h$. We introduce the
non-dimensional constant $\nu = \mu T / L^2$, which is a measurement
of the quotient of viscous by inertial forces in our flow, and
measures in a physically appropriate manner how far from ideal our
flow really is. The non-dimensionalization procedure simply means
introducing new variables

\[\widetilde{u}^1(y,s) \equiv \frac{T}{L} u^1(Ly_1,hy_2,Ts), \]
\[\widetilde{u}^2(y,s) \equiv \frac{T}{h} u^2(Ly_1,hy_2,Ts),\]
and
\[\widetilde{p}(y,s) \equiv \frac{T^2}{L^2} p(Ly_1,hy_2,Ts).\]

which results in the system

\begin{equation} \label{adim}
\left\{
\begin{array}{l}
\partial_s \widetilde{u}^1 + \widetilde{u} \cdot \nabla_y \widetilde{u}^1 = - \partial_{y_1} \widetilde{p} +
\nu\left[ \partial^2_{y_1}\widetilde{u}^1 + \displaystyle{\frac{L^2}{h^2}} \partial^2_{y_2}\widetilde{u}^1\right]\\ \\
\partial_s \widetilde{u}^2 + \widetilde{u} \cdot \nabla_y \widetilde{u}^2 =
- \displaystyle{\frac{L^2}{h^2}}\partial_{y_2} \widetilde{p} +
\nu\left[ \partial^2_{y_1}\widetilde{u}^2 +
\displaystyle{\frac{L^2}{h}} \partial^2_{y_2}\widetilde{u}^2\right]\\ \\
\mbox{div}_y \widetilde{u} = 0.
\end{array}
\right.
\end{equation}

We introduce $\vare \equiv h/L$, which is assumed to be a small,
non-dimensional, parameter because we are focusing in a thin layer.
We also want to consider $\nu$ small. A key issue in boundary layer
theory is that the magnitude of the small parameters $\vare$ and
$\nu$ are naturally related. Indeed, if $\nu \ll \vare^2$ then
matched asymptotics indicates that, to leading order,
$\widetilde{u}$ satisfies the Euler system (with pressure
independent of the vertical variable) and with no-slip conditions at
the boundary. These boundary conditions are inconsistent with the
fact that the Euler system is of first order. On the other hand, if
$\nu \gg \vare^2$ then, to leading order, $\widetilde{u}$ is such
that $\partial^2_{y_2} \widetilde{u} = 0$, which, together with the
no-slip boundary condition implies that $\widetilde{u} = c(y_1)y_2$,
for some vector $c$. Now, if $\widetilde{u}$ is to represent the
behavior of the flow in a thin layer near the boundary, then the
velocity $\widetilde{u}$ should match the inviscid velocity in the
limit $y_2 \to \infty$, and not just blow-up. The only regime that
appears to yield a consistent asymptotic regime is
\begin{equation} \label{blthcknss} \nu/\vare^2 = \mathcal{O}(1).
\end{equation}
From another perspective, condition \eqref{blthcknss} highlights the
region near the boundary where the vertical viscous stress balances
the inertial terms in the Navier-Stokes system. Assuming $\nu =
\vare^2$ and implementing matching asymptotics for $\nu$ small we
obtain the following system for the leading order approximation,
denoted $v = (v^1,v^2)$,
\begin{equation} \label{prandtl}
\left\{ \begin{array}{l}
\partial_s v^1 + v \cdot \nabla_y v^1 = - \partial_{y_1} q +
\partial^2_{y_2}v^1\\
\partial_{y_2} q =0 \\
\mbox{div}_y v = 0.
\end{array} \right.
\end{equation}

These are the {\it unsteady Prandtl equations} for the boundary
layer profile $v$. They represent the behavior of the flow near the
boundary. To obtain a complete problem, these equations must be
supplemented with boundary conditions. First we impose the no-slip
condition
\[v = 0 \mbox{ at } y_2 = 0.\]
An additional condition must be imposed in order to capture the
assumption that far from the boundary layer the small viscosity
Navier-Stokes solutions match the Euler solutions. Let $u^{\nu}$ a
family of solutions of the non-dimensional Navier-Stokes equations
with non-dimensional viscosity $\nu$ and $u^E$ be a solution of the
incompressible Euler equations. For example, can assume that both
the family $u^{\nu}$ and $u^E$ are defined by solving the
Navier-Stokes equations and the Euler equations with the same
initial data $u^{\nu}(x,0) = u^E(x,0) = u_0(x)$.

Going back to the Prandtl system, we expect that $v(y_1,y_2,t) \to
u^E(y_1,0,t) \equiv U(y_1,t)$, as $y_2 \to \infty$. Let $p^E$ be the
pressure associated with the Euler solution $u^E$.  Since $q$ does
not depend on $y_2$, looking at $y_2 \to \infty$ makes it also
natural to assume that $q(y_1,t) = p^E(y_1,0,t)$. If we look at the
Euler equations and evaluate them at $x_2 = 0$ we obtain the
following relation:
\[U_t + U U_{y_1} = - P^E_{y_1} = -q_{y_1}, \]
which is called Bernoulli's Law. This means that, for the Prandtl
equation \eqref{prandtl}, the condition at infinity $U$ determines,
up to an irrelevant constant, the pressure $q$.

With this construction, we hope that, when $\nu$ is small,
\begin{equation} \label{bdlayapprox}
u^{\nu}(x_1,x_2,t)  = \left\{ \begin{array}{l}
(v^1,\vare v^2)(x_1,x_2/\vare,t) \mbox{ for } x_2 < \lambda(\nu) \\
u^E(x_1,x_2,t) \mbox{ for } x_2 > \lambda(\nu) \end{array} \right. +
o(1),
\end{equation}
where $\lambda(\nu)$ is any cutoff distance such that $\vare \ll
\lambda(\nu) \ll 1$.

In addition to the time-dependent Prandtl equation, this derivation
also yield the {\it steady Prandtl equation}, given by
\begin{equation} \label{sprandtl}
\left\{ \begin{array}{l}
v \cdot \nabla_y v^1 = - \partial_{y_1} q +
\partial^2_{y_2}v^1\\
\partial_{y_2} q =0 \\
\mbox{div}_y v = 0.
\end{array} \right.
\end{equation}
The typical problem associated with this equation is a quarter-plane
BVP, where a profile $v(0,y_2)$ is given and one attempts to find
$v(y_1,y_2)$ for $y_1>0$ as the induced boundary layer profile over
a half-plane plate.

The derivation above is a nice example of multiscale asymptotic
analysis, and from the complicated issue surrounding the interaction
of nearly inviscid flow with material boundaries it derives a new
equation, \eqref{prandtl}, and a simplified asymptotic model for the
behavior of Navier-Stokes solutions near a material boundary, given
by \eqref{bdlayapprox}. The key issue that such a model raises is
its validity (mathematical) and applicability (physical).

This model has been found useful in applications, specially where it
concerns laminar boundary layers, and when its usefulness begins to
break down, suitable extensions of the model have been obtained,
notably the so called ``triple deck'' expansions, where an
intermediate thin layer is added between the viscosity-dominated
internal layer and the free irrotational Euler flow. In this
intermediate layer, the flow is ideal, but not necessarily
irrotational. One important situation where the boundary layer
ansatz breaks down is when {\it boundary layer separation} occurs.
Recall that one of the assumptions in deriving the Prandtl equation
was that the disparity between ideal and viscous flow be
concentrates in a thin layer {\it near the boundary}. It is quite
common, even well within laminar flow regimes, that the boundary
layer detaches itself from the boundary and affects the bulk of the
flow. In that case, Prandtl`s theory and its extensions break down
as models. In the next section we will consider what is known
regarding the rigorous validation of the asymptotic approximation
described here.

\section{Prandtl`s equation}

The purpose of this Section is to present the known theory for
Prandtl's equation, without much detail. The first question one must
address regarding an approximate model is: can I solve it? The
initial and boundary value problems for the Euler and Navier-Stokes
equations are well-posed, globally in the case of the half-plane
with reasonable initial data. So the issue of whether Euler $+$
Prandtl is a good approximation for Navier-Stokes with $\nu$ small,
in the sense discussed in the previous section, depends first on
understanding the well-posedness for Prandtl`s equation. The
mathematical theory of the Prandtl equation only got started in the
sixties, by O. Oleinik. Over the years, Oleinik and her group made a
large number of contributions to the theory of Prandtl`s equation
and its many variants, collected and explained in the book
\cite{oleinikbook}. For the present discussion, we would like to
focus on one specific result, that first appeared in \cite{oleinik}.
We also refer the reader to the survey \cite{E00}, for the
discussion of Oleinik`s result and its relation with the blow-up
result of W. E and B. Engquist, \cite{EE97}.

Let $v = (v^1,v^2)$ be a solution of the IBVP for Prandtl`s
equation, which we write as

\begin{equation} \label{IBVPprandtl}
\left\{ \begin{array}{l}
\partial_t v^1 + v \cdot \nabla v^1 = -\partial_{x_1}q +
\partial^2_{x_2}v^1\\
\mbox{div} v = 0\\
v(x_1,0,t) = 0 \mbox{ and } \lim_{x_2 \to \infty} v^1(x,t) =
U(x_1,t)\\
v(x,0) = v_0(x),
\end{array} \right.
\end{equation}

Where $-\partial_{x_1}q = U_t + UU_{x_1}$. The result we wish to
discuss is the following

\begin{theorem} \label{oleinik} (Oleinik 1967) Assume that both
$U$ and $v_0^1$ are positive and that, in addition, $\partial_{x_2}
v_0^1 \geq 0$. Then there exists a unique global strong solution of
\eqref{IBVPprandtl}.
\end{theorem}

We will not present a proof of this result, but we will discuss a
key part of the proof, which is the recasting of this problem as a
scalar, degenerate parabolic scalar equation, using Crocco`s
transformation.

We begin by taking the derivative of Prandtl`s equation with respect
to $x_2$ and introduce the new dependent variable $\omega(x,t)
\equiv \partial_{x_2} v^1$. System \eqref{IBVPprandtl} is equivalent
to the following IBVP:

\begin{equation} \label{vortprandtl}
\left\{ \begin{array}{l}
\partial_t \omega + v \cdot \nabla \omega = \partial^2_{x_2}\omega\\
v = K[\omega] \\
(\partial_{x_2}\omega)(x_1,0,t) = \partial_{x_1}q  \mbox{ and } \lim_{x_2 \to \infty} \omega(x,t) = 0 \\\
\omega(x,0) = \partial_{x_2}v_0(x),
\end{array} \right.
\end{equation}
where the vector operator $K$ reconstructs $(v^1,v^2)$ by first
integrating in the vertical variable to obtain $v_1$ and then using
the divergence-free condition and integrating again in the vertical
variable to obtain $v_2$. The new equation \eqref{vortprandtl} is,
in a sense, the vorticity formulation of problem
\eqref{IBVPprandtl}.

We assume that the solution $v = v(x,t)$ we seek satisfies the
condition $\partial_{x_2} v^1(x,t) > 0$ for all $x \in \mathbb{H},
t>0$. In particular, this means that, for each fixed $(x_1,t)$, the
map $x_2 \mapsto v^1(x_1,x_2,t)$ is invertible. Let us denote this
inverse by $h = h(x_1,\xi,t)$. In other words, we have
\begin{equation} \label{hdef}
v^1(x_1,h(x_1,\xi,t),t) = \xi, \mbox{ for all } \xi > 0.
\end{equation}

The Crocco`s transform consists of introducing the new dependent
variable:
\begin{equation} \label{Wdef}
W = W(x_1,\xi,t) \equiv \omega(x_1,h(x_1,\xi,t),t),
\end{equation}

We verify that $W$ is a solution of the following IBVP:

\begin{equation} \label{croccoprandtl}
\left\{ \begin{array}{l}
\partial_t W + \xi W_{x_1} - (\partial_{x_1} q) W_{\xi} = W^2 \partial^2_{\xi}W\\
W W_{\xi} = q_{x_1}  \mbox{ for } \xi = 0 \\
W(x,0) = (\partial_{x_2}v_0)(x_1,h(x_1,\xi,0))
\end{array} \right.
\end{equation}

Indeed, we can compute directly to obtain:

\[ \partial_t W = -\frac{v^1_t}{\omega} \omega_{x_2} + \omega_t
\mbox{; } \partial_{x_1} W = -\frac{v^1_{x_1}}{\omega} \omega_{x_2}
+ \omega_{x_1} \mbox{; }
\partial_{\xi} W = \frac{\omega_{x_2}}{\omega} \mbox{; }
\partial^2_{\xi} W = \frac{\omega_{x_2x_2}}{\omega^2} -
\frac{\omega_{x_2}^2}{\omega^3}.\]

Substituting the corresponding equalities above into
\eqref{croccoprandtl}, and using \eqref{Wdef}, \eqref{hdef} and the
evolution equations in \eqref{IBVPprandtl} and \eqref{vortprandtl}
verifies the evolution equation in \eqref{croccoprandtl}. In
addition, we can check directly that
\[\omega_{x_2}(x_1,x_2,t) = W(x_1, g^{-1}(x_1,x_2,t),t)
W_{\xi}(x_1,g^{-1}(x_1,x_2,t),t),\] which, together with the
boundary condition in \eqref{vortprandtl} gives the boundary
condition in \eqref{croccoprandtl}. Problem \eqref{croccoprandtl} is
a scalar, degenerate parabolic equation, which is amenable to fairly
standard treatment, using fixed point methods, and satisfies a
comparison principle. In particular, the sign of $W$ is retained in
the evolution, and therefore, the monotonicity condition on $v^1$,
necessary for the validity of Crocco`s transform, is retained as
well.

Finally, once a solution $W$ is obtained for problem
\eqref{croccoprandtl}, one must reconstruct a solution to the
original problem. Recall that
\[v^1(x,g(x_1,\xi,t),t) = \xi,\]
and therefore, differentiating this identity with respect to $\xi$
gives
\[ \omega(x_1,g(x_1,\xi,t),t) \frac{\partial g}{\partial
\xi}(x_1,\xi,t) = 1.\] Recalling \eqref{Wdef}, this implies that
\[ g(x_1,\xi,t) = \int_0^{\xi} (W(x_1,\eta,t))^{-1} d\eta, \]
which allows the reconstruction of $\omega$ from $W$ by means of
\eqref{Wdef}. The interested reader may prove, as an exercise, that
such an $\omega$ is, in fact, a solution of \eqref{vortprandtl}.

This result, and others proved by Oleinik and her group, give
useful, rigorously established descriptions of the vanishing
viscosity asymptotics, but depend, to a greater or lesser extent, on
monotonicity conditions such as $\partial_{x_2} v_0^1 \geq 0$. As we
have seen, the monotonicity assumption is needed for the validity of
the Crocco`s transformation, but this assumption might just be a
feature of the method, rather than an essential limitation of the
theory. In 1997, E and Engquist produced a counterexample which
showed that Prandtl`s equation develops finite-time singularities if
the monotonicity condition is not imposed, see \cite{EE97}. In fact,
E and Engquist`s example suggests that the role of the monotonicity
assumption is to prevent boundary layer separation, a phenomenon
that actually occurs in real flows and corresponds to a breakdown of
the Prandtl ansatz.

An alternative to the half-plane analysis described above is to
study well-posedness of Prandtl`s equation in bounded intervals in
$x_1$, where the horizontal velocity of the boundary layer profile
is specified in one side of the interval, and the length of the
interval or the time horizon of the analysis are chosen small enough
to prevent separation. Such a result was first proved by Oleinik in
\cite{oleinik66}. Recently, Z. Xin and L. Zhang improved Oleinik`s
result, showing existence of a global (in time) weak solution for
Prandtl`s equation on a finite horizontal interval, if the pressure
is {\it favorable}, i.e., $q_{x_1} < 0$, see \cite{XZ04} and
\cite{zz06}. This condition is also known to discourage boundary
layer separation.

A different approach to the theory of Prandtl`s equation was taken,
initially by A. Asano, in a couple of unpublished manuscripts, and
later by R. Caflisch and M. Sammartino, in a pair of articles, see
\cite{SC98}, recently further improved by Lombardo, Cannone and
Sammartino in \cite{LCS03}. The basic idea is that, without the
monotonicity condition, or something analogous to it, one expects
the initial-boundary value problem \eqref{IBVPprandtl} to be
ill-posed. As a result, it becomes natural to look for local (in
time) solutions for Prandtl`s equation in analytic function spaces,
using results of Cauchy-Kowalewska type. The main results in
\cite{SC98} were well-posedness of the problem \eqref{IBVPprandtl}
if the data $v_0$ and $U$ are analytic, and compatible. In
\cite{LCS03} the analyticity requirement on $v_0$ was imposed only
in the horizontal variables.

Of course, the well-posedness in analytic spaces, and the blow-up
example by E and Engquist does not prove that \eqref{IBVPprandtl} is
ill-posed, which at this time remains an interesting open problem.

To conclude this section, it would make sense to mention the
contribution of E. Grenier, which he describes as nonlinear
instability of the Prandtl boundary layer. His result is not about
Prandtl's equation per se, but about the vanishing viscosity limit
of the Navier-Stokes equations. His result can be interpreted as
mathematical evidence that the Prandtl ansatz is not always valid
for solutions of the Navier-Stokes system in the half-plane with
small viscosity, see \cite{grenier00}. In other words, although the
theory of Prandtl`s equation is relevant for understanding the
vanishing viscosity limit, there is more to the original observation
(O) than Prandtl`s original explanation for it.

\section{Kato`s condition}

In this section we move away from Prandtl`s equation, and we begin a
study of the vanishing viscosity limit from a broader point of view.
Our first observation should be that, even in the absence of
boundaries, all the mysteries of turbulence lurk in the background
of the vanishing viscosity limit, see for example \cite{LMN06} for a
small part of this story. However, under moderate regularity
assumptions, for example, if the initial vorticity is bounded,
explicit estimates for the difference between Euler and small
viscosity Navier-Stokes solutions are known, see \cite{chemin96}.
Also, and this distinction is a key point here, for very irregular
flow, we still have the existence of subsequences of solutions of
small viscosity Navier-Stokes converging to weak solutions of the
Euler equations, up to and including initial vorticities which are
measures, see \cite{Majda93,LNX01}, but then no estimates on the
difference are provided, or expected. Basically, in the absence of
boundaries, as long as the underlying ideal flow has enough
regularity so that uniqueness of weak solutions to the Euler
equations is known, we have actual convergence of the vanishing
viscosity limit. Furthermore, as long as existence of weak solutions
is known we also have compactness of the vanishing viscosity
sequence and weak continuity of the Euler/Navier-Stokes
nonlinearity. Nonuniqueness of weak solutions for Euler equations is
also known, see the remarkable paper \cite{dellelis}, and references
therein, for the current knowledge on this nonuniqueness, but the
behavior of the vanishing viscosity limit for these examples is a
very interesting open problem.

As soon as we consider flows in the presence of boundaries, the
story changes quite dramatically. Very little is actually known
mathematically, and this very little is precisely the object under
discussion in these notes. Physically, boundaries are the most
natural source of small scales in incompressible flows, precisely
through the boundary layer mechanism, and these small scales are the
source of the irregularities that justify considering irregular 2D
flows in the first place. The point of departure in our discussion
will be a classical open problem, which we formulate below.

\vspace{0.5cm}

{\bf $\partial$ Layer Problem:} Let $u^{\nu}$ be a sequence of
solutions of the incompressible Navier-Stokes equations in two space
dimensions, in a smooth bounded domain $\Omega$, satisfying the
no-slip boundary condition on $\partial \Omega$, with initial data
$u_0^{\nu}$, bounded in $L^2$. Is there a subsequence $u^{\nu_k}$
converging weakly in $L^2$ to a vector field $u$, which is a weak
solution of the incompressible Euler equations in $\Omega$ with some
initial data $u_0 = \lim u_0^{\nu_k}$?

\vspace{0.5cm}

This problem is open even if $\omega_0^{\nu} = \omega_0 \in
C_c^{\infty}(\Omega)$, with $\omega_0 = \mbox{ curl }u_0$. Let us
focus, for simplicity, in this case. Clearly, the Navier-Stokes
equations have a unique smooth solution $u^{\nu}$ with initial data
$u_0$, and the Euler equations also have a unique smooth solution
$u$ with the same initial data. We will see that there are examples
where $u^{\nu} \to u$ in $L^2$, but the answer to the problem above
may be positive even when $u^{\nu}$ does not converge to $u$,
because there may be weak solutions of the incompressible Euler
equations with initial velocity $u_0$ which are not $u$.

In 1984, T. Kato wrote a short note where he proved a sharp
criterion for the convergence of $u^{\nu}$ to $u$, see
\cite{Kato84}. The observation by Kato is remarkable for at least
two reasons. First, as we shall see, it is very natural from the
analytical point of view. Second, it places the condition for
convergence on the behavior of the small viscosity sequence at a
distance $\mathcal{O}(\nu)$ of the boundary of $\Omega$, hence in a
much smaller region than what is the natural domain of the boundary
layer. Next, we state and prove a simple version of Kato`s
criterion.

We focus in the case $\Omega = \{|x| < 1\}$ in $\real^2$. Let
$\omega_0 \in C^{\infty}_c(\Omega)$ and $u_0 \equiv K[\omega_0]$,
where $K$ is the Biot-Savart operator in the unit disk. Let
$u^{\nu}$ be the unique classical solution of the Navier-Stokes
equation in $\Omega$ with no-slip boundary condition and initial
velocity $u_0$, and $u$ be the unique smooth solution of the Euler
equations with $u \cdot x = 0$ for $|x|=1$ and initial velocity
$u_0$.

\begin{theorem} (Kato 1984) Fix $T>0$. There exists a constant $c>0$ such that
$u^{\nu} \to u$ strongly in $L^{\infty}((0,T);L^2(\Omega))$ if and
only if $\nu \int_0^T \| \nabla u^{\nu}(\cdot,t)
\|_{L^2(\Gamma_{c\nu})}^2 dt \to 0$ as $\nu \to 0$, where
$\Gamma_{c\nu} \equiv \{1-c\nu < |x| < 1\}$.
\end{theorem}

\begin{proof}
First consider the energy identities for both $u^{\nu}$ and $u$. We
have, for each $t>0$,
\[ \|u^{\nu}(\cdot,t)\|^2_{L^2(\Omega)} = \|u_0\|^2_{L^2(\Omega)}
+ \nu \int_0^t \int_{\Omega} |\nabla u^{\nu}|^2 dxdt, \] and
\[ \|u(\cdot,t)\|^2_{L^2(\Omega)} = \|u_0\|^2_{L^2(\Omega)}. \]
Therefore, if $u^{\nu} \to u$ strongly in
$L^{\infty}((0,T);L^2(\Omega))$, then
$\|u^{\nu}(\cdot,t)\|^2_{L^2(\Omega)} \to
\|u(\cdot,t)\|^2_{L^2(\Omega)}$ for almost all time, and therefore
\[\nu \int_0^t \int_{\Omega} |\nabla u^{\nu}|^2 dxdt \to 0,\]
for each $t>0$, not almost everywhere anymore since the integral in
time is increasing in time, and therefore,
\[\nu \int_0^t \int_{\Gamma_{c\nu}} |\nabla u^{\nu}|^2 dxdt \to 0,\]
as we wished.

To prove the other implication, fix $\vare > 0$ and let
$\phi^{\vare} \in C^{\infty}_c(\Omega)$ be such that
$\phi^{\vare}(x) = \varphi^{\vare}(|x|)$, with $\varphi^{\vare}(x) =
1$ for $|x| < 1-\vare$, $\varphi^{\vare}(x) = 0$ for $1 - \vare/2 <
|x| \leq 1$, and $\varphi^{\vare}$ decreases monotonically from $1$
to $0$. Let $\omega = \mbox{ curl }u$ be the vorticity and $\psi$ be
the stream function associated with the Euler flow $u$. Define
 \[u_{\vare} = \nabla^{\perp}(\phi^{\vare} \psi) = (-\partial_{x_2}(\phi^{\vare}
\psi),\partial_{x_1}(\phi^{\vare} \psi)).\]

Let $v_{\vare} \equiv u - u_{\vare} =
\nabla^{\perp}((1-\phi^{\vare})\psi)$. The stream function $\psi$
vanishes at $|x|=1$, and it can be assumed to be uniformly bounded
in $C^k$, for any $k = 1,2, \ldots$, so that we can easily obtain
the following estimates on $v_{\vare}$:

\begin{equation} \label{ev1} \|v_{\vare}\|_{L^{\infty}((0,T);L^2(\Omega))} \leq C
\vare^{1/2} \end{equation}

\begin{equation} \label{ev2}\|\partial_t v_{\vare}\|_{L^1((0,T);L^2(\Omega))} \leq C
\vare^{1/2}\end{equation}

\begin{equation} \label{ev3}\|\nabla v_{\vare}\|_{L^{\infty}((0,T);L^2(\Omega))} \leq C
\vare^{-1/2}\end{equation}

\begin{equation} \label{ev4}\|v_{\vare}\|_{L^{\infty}((0,T) \times \Omega)} \leq C\end{equation}

\begin{equation} \label{ev5}\|\nabla v_{\vare}\|_{L^{\infty}((0,T) \times \Omega)} \leq C
\vare^{-1}\end{equation}

In addition, we require certain estimates on $u^{\nu}$, uniform in
$\nu$, which we collect below

\begin{equation} \label{eunu1} \|u^{\nu}\|_{L^{\infty}((0,T);L^2(\Omega))} \leq
C \end{equation}

\begin{equation} \label{eunu2} \nu \|\nabla u^{\nu}\|^2_{L^{\infty}((0,T);L^2(\Omega))} \leq
C. \end{equation}

By using Cauchy-Schwarz in time, we also have

\begin{equation} \label{eunu3} \nu^{1/2} \|\nabla
u^{\nu}\|_{L^1((0,T);L^2(\Omega))} \leq C T^{1/2} \left( \nu
\int_0^T \|\nabla u^{\nu} \|_{L^2(\Omega)}^2 dt \right)^{1/2} \leq C
\end{equation}

Finally, a version of Poincar\'{e}'s Inequality, which reads

\begin{equation} \label{poincare}
\|u^{\nu}\|_{L^2(\Gamma_{\vare})} \leq C \vare\|\nabla u^{\nu}
\|_{L^2(\Gamma_{\vare})}
\end{equation}

Now we estimate, omitting the explicit dependence of $v$ on $\vare$:

\[\|u^{\nu} - u \|_{L^2(\Omega)}^2 =  \|u^{\nu}\|_{L^2(\Omega)}^2 + \|u
\|_{L^2(\Omega)}^2 - 2 \langle u^{\nu},u \rangle \]
\[\leq 2 \|u_0 \|_{L^2(\Omega)}^2 - 2\langle u^{\nu},u-v \rangle - 2\langle u^{\nu},v
\rangle.\]

We have that
\[|\langle u^{\nu},v \rangle| \leq
\|u^{\nu}\|_{L^{\infty}((0,T);L^2(\Omega))}
\|v\|_{L^{\infty}((0,T);L^2(\Omega))} \leq C \vare^{1/2}.\]

And therefore, taking $\vare = C\nu$,

\begin{equation} \label{stp1} \|u^{\nu} - u \|_{L^2(\Omega)}^2\leq
2 \|u_0 \|_{L^2(\Omega)}^2 - 2\langle u^{\nu},u-v \rangle +
\mathcal{o}(1).
\end{equation}

We multiply the Navier-Stokes equation by $u-v$, integrate in space
and time, and integrate by parts to obtain
\[ \langle u^{\nu},u-v \rangle - \langle u_0,u_0-v \rangle =
\int_0^t \langle u^{\nu}, u^{\nu} \cdot \nabla (u-v) \rangle - \nu
\langle \nabla u^{\nu},\nabla(u-v) \rangle + \langle
u^{\nu},\partial_t(u-v) \rangle dt, \]

which, together with \eqref{stp1} implies

\begin{equation} \label{stp2} \|u^{\nu} - u \|_{L^2(\Omega)}^2\leq
\int_0^t [-2 \langle u^{\nu}, u^{\nu} \cdot \nabla (u-v) \rangle +2
\nu \langle \nabla u^{\nu},\nabla(u-v) \rangle -2 \langle
u^{\nu},\partial_t(u-v) \rangle ] dt + o(1).
\end{equation}

We write
\[ \int_0^t \langle u^{\nu},\partial_t(u-v) \rangle  dt = -\int_0^t
\langle u^{\nu}, u \cdot \nabla u \rangle dt + \int_0^t  \langle
u^{\nu}, \partial_t v \rangle dt,\]
and we have
\[\left| \int_0^t \langle u^{\nu}, \partial_t v \rangle dt \right|
\leq \int_0^t \|u\|_{L^2} \|\partial_t v\|_{L^2} \leq C \nu^{1/2},\]
hence,

\begin{equation} \label{stp3}\|u^{\nu} - u \|_{L^2(\Omega)}^2\leq
\int_0^t [-2 \langle u^{\nu}, u^{\nu} \cdot \nabla (u-v) \rangle +2
\nu \langle \nabla u^{\nu},\nabla(u-v) \rangle + 2 \langle u^{\nu} ,
u \cdot \nabla u\rangle] dt + o(1).
\end{equation}

Note that
\[ \langle (u^{\nu} - u) , (u^{\nu} - u) \cdot \nabla u\rangle
= \langle u^{\nu} , u^{\nu} \cdot \nabla u \rangle - \langle u^{\nu}
, u \cdot \nabla u\rangle ,\]

and therefore, \eqref{stp3} becomes

\begin{equation} \label{stp4}\|u^{\nu} - u \|_{L^2(\Omega)}\leq
2  \int_0^t[-\langle (u^{\nu}-u), (u^{\nu} - u) \cdot \nabla u
\rangle + \nu \langle \nabla u^{\nu},\nabla(u-v) \rangle + \langle
u^{\nu} , u^{\nu} \cdot \nabla v \rangle] dt + o(1).
\end{equation}

We analyze each term in the right hand side of \eqref{stp4}
separately. First,

\[ \left|\int_0^t \langle (u^{\nu}-u), (u^{\nu} - u) \cdot \nabla u
\rangle dt \right| \leq \int_0^t \|u^{\nu} - u\|^2_{L^2(\Omega)}
|\nabla u| dt \leq C \int_0^t \|u^{\nu} - u\|^2_{L^2(\Omega)}dt.\]

Second,

\[\left| \nu \int_0^t \langle \nabla u^{\nu},\nabla(u-v) \rangle dt
\right| \leq \left| \nu \int_0^t \langle \nabla u^{\nu},\nabla u
\rangle dt \right| + \left| \nu \int_0^t \langle \nabla
u^{\nu},\nabla v \rangle dt \right| \]
\[\leq \nu \int_0^t \|\nabla u^{\nu}\|_{L^2(\Omega)} \|\nabla u\|_{L^2(\Omega)}
dt  + \nu \int_0^t\|\nabla u^{\nu}\|_{L^2(\Gamma_{c\nu})} \|\nabla
v\|_{L^2(\Gamma_{c\nu})}\] \[ \leq C \nu \|\nabla
u^{\nu}\|_{L^1((0,T);L^2(\Omega))} +  C \nu^{1/2} \|\nabla u^{\nu}
\|_{L^1((0,T);L^2(\Gamma_{c\nu}))} \]
\[ \leq C \nu^{1/2} t^{1/2} \|\nabla u^{\nu}
\|_{L^2((0,T);L^2(\Gamma_{c\nu}))} + o(1) = Ct^{1/2} \left( \nu
\int_0^t \|\nabla u^{\nu} \|_{L^2(\Gamma_{c\nu})}^2 \right)^{1/2} =
o(1), \] by Kato's criterion.

Third,
\[\left| \int_0^t \langle u^{\nu},u^{\nu}\cdot \nabla v \rangle dt
\right| = \left| \int_0^t \langle v,u^{\nu}\cdot \nabla u^{\nu}
\rangle dt \right| \]
\[ \leq 2 \|v\|_{L^{\infty}((0,T) \times \Omega)} \int_0^t
\|u^{\nu}\|_{L^2(\Gamma_{c\nu})} \|\nabla
u^{\nu}\|_{L^2(\Gamma_{c\nu})} dt \leq C \nu \int_0^t \|\nabla
u^{\nu}\|_{L^2(\Gamma_{c\nu})}^2 dt = o(1), \] where we used
Poincar\'{e}'s Inequality and Kato's criterion in the last line.

Therefore, we get back to \eqref{stp4} and use the estimates
obtained above to conclude that

\[\|u^{\nu} - u \|_{L^2(\Omega)}^2\leq C \int_0^t \|u^{\nu} -
u\|^2_{L^2(\Omega)}dt + o(1), \] which, by Gronwall, concludes the
proof.

\end{proof}

A number of observations are in order. First, the proof above can be
adapted to general bounded domains in $\real^n$, to Leray solutions
of the Navier-Stokes equations and to problems with forcing and
different initial data, as long as the forcing and the initial data
converge to the corresponding ones when viscosity vanishes. The two
equivalent conditions in the statement of the Theorem are also
equivalent to weak convergence of $u^{\nu}$ to $u$ in $L^2$,
pointwise in time. Clearly the strong convergence implies the weak
convergence, and we weak convergence implies that
\[\nu \int_0^t \|\nabla u^{\nu}\|_{L^2(\Omega)}^2 dt \to 0, \mbox{
as } \nu \to 0,\] which can be seen using the energy identities and
the weak lower semicontinuity of the $L^2$ norm. All these
observations were contained in the original work by Kato,
\cite{Kato84}.

There have been several results modifying, extending and improving
Kato's criterion in the literature: work by R. Temam and Xiaoming
Wang, \cite{TW97,TW98}, where Kato`s condition is replaced by a
condition on integrability of pressure near the boundary, work by
Xiaoming Wang, \cite{Wang01} where the criterion is cast in terms of
integrability of derivatives of tangential components of velocity,
and work by James Kelliher, \cite{kelliher07}, where the
integrability condition is placed on vorticity.

As we will see in the last lecture, there are some (trivial)
examples where Kato`s criterion holds true, but this is not known in
general, and it would not be surprising if it turns out to be false
most of the time. Note that, if Kato`s criterion fails, this implies
non-vanishing energy dissipation occurring in a thin neighborhood
around the boundary, as viscosity vanishes.

Given the disparity of the understanding of the vanishing viscosity
limit in the presence or absence of boundaries, it becomes natural
to ask how much boundary it is required to obstruct the proof of
convergence of the vanishing viscosity limit. The author has been
involved with two recent results in this direction. First, with D.
Iftimie and H. Nussenzveig Lopes, we have proved that if
$u^{\nu,\vare}$ is the solution of the Navier-Stokes equations in
the exterior of the domain $\vare \Omega$, with no-slip boundary
conditions and fixed initial vorticity and $u$ is the solution of
the Euler equations in the full space, with the same initial
vorticity then there exists a constant $c>0$ such that, if $\vare <
c\nu$, then $u^{\nu,\vare} \to u$ strongly in $L^2$, see
\cite{ILN08}. Second, with J. Kelliher and H. Nussenzveig Lopes, we
have proved that if $u^{\nu,R}$ is the solution of the Navier-Stokes
equations on a domain $\Omega_R$ which contains the ball of radius
$R$, with fixed initial vorticity, and again $u$ is the solution of
the Euler equation with the same initial vorticity. We have proved
that, as long $R \to \infty$ as $\nu \to 0$, then $u^{\nu,R}$
converges strongly to $u$ as $\nu \to 0$, see \cite{KLN08}. Both
proofs follow the basic idea and the structure of the proof of
Kato's criterion as we presented it here.

\section{Vanishing viscosity with Navier friction $\partial$ condition}

The main purpose of this Section is to highlight the role of
vorticity production by the interaction of the flow with the
boundary in the difficulty surrounding the Boundary Layer Problem.
In particular, we briefly discuss the vanishing viscosity limit for
flows in the full space, and also for flows in a bounded domain
satisfying Navier friction condition at the boundary. Let us begin
with the full space case.

A basic ingredient for the proof of compactness of the vanishing
viscosity limit in the style of \cite{Majda93} is the set of
estimates obtained through the vorticity formulation of the
Navier-Stokes equations. Vorticity play a fundamental role in the
analysis of incompressible fluid flow. This is specially true in two
dimensions, because vorticity is conserved along particle
trajectories for ideal flow, and it satisfies a nice
convection-diffusion equation for viscous flow, namely:
\begin{equation} \label{vorteq}
\left\{ \begin{array}{l} \omega_t + u \cdot \nabla \omega =
\nu\Delta \omega\\
u = K[\omega],
\end{array} \right.
\end{equation}
where $K$ is the Biot-Savart law associated with the flow domain. In
the absence of material boundaries, for example, in the case of the
full plane, equation \eqref{vorteq} implies a priori $L^p$ estimates
for vorticity for any $1 \leq p \leq \infty$, as long as the initial
vorticity is in $L^p$. The Biot-Savart operator is a
pseudo-differential operator of order $-1$, which implies
compactness of the sequence of velocities in $L^q$, for any $q <
p^{\ast}$, $p^{\ast} = 2p/(2-p)$ (assuming $1 \leq p < 2$), and the
nonlinearity of the Euler equations being quadratic, all that is
required is compactness in $L^2$ in order to prove that the flows
obtained through the limit of vanishing viscosity satisfy Euler
equation. The case $p=1$, which is what is involved in the vortex
sheet initial data problem, is critical for the argument above, and
requires a more involved analysis.

To highlight the importance of the added compactness provided by
vorticity estimates, we consider briefly the situation without it.
In that case, the only a priori estimate available comes from the
energy estimate, which gives solely an $L^2$-bound on velocity,
independent of viscosity. By Banach-Alaoglu, an uniform $L^2$-bound
can provide us with a sequence of approximations which converge {\it
weakly} in $L^2$ to a "candidate" for limit solution $u$. To prove
that $u$ really is a solution of Euler's equation, one must pass to
the limit in the nonlinearity $u^n \cdot \nabla u^n$, which, after
integration by parts on a weak formulation, can be written in terms
of {\it quadratic} expressions in the components of velocity. Now,
the point is: quadratic expressions are not continuous with respect
to the weak topology of $L^2$. Two classical examples (in one space
dimension) are:
\[ f^n = \left\{ \begin{array}{l} \frac{1}{\sqrt{n}}, \mbox{ for }
|x| \leq \frac{1}{n} \\ 0 \mbox{ for } |x| > \frac{1}{n} \end{array}
\right. \mbox{ \it (concentration)},\] and
\[ g^n = sin(nx) \mbox{ \it (oscillation)}. \]
Both these sequences converge weakly to zero in $L^2$, but their
squares converge to something else - we leave the details as an
exercise to the reader. The bottom line is this: strong convergence
in $L^2$ is required to pass to the limit in quadratic expressions,
and this strong convergence has to come from some additional a
priori estimate.

The presence of vorticity estimates breaks down in the presence of
material boundaries because, although vorticity still satisfies
equation \eqref{vorteq} in a domain with boundary, there is no
natural boundary condition for vorticity to complete equation
\eqref{vorteq}, beyond the nonlocal boundary condition $K[\omega]
\cdot \widehat{\tau} = 0$ (only the vanishing of the tangential
component of velocity is required to vanish because the vanishing of
the normal component of velocity is implicit in the Biot-Savart
law). In particular, one may argue that it is {\em the lack of
control on the production of vorticity at the boundary} which
obstructs the solution of the $\partial$ Layer Problem.

Mathematicians are not discouraged by difficulties - if we cannot
solve a problem, we find another problem nearby which we can
actually solve. In the book \cite{lions69}, J.-L. Lions used as an
illustration the following problem:

\begin{equation} \label{fbnseq}
\left\{  \begin{array}{l}
\partial_t u + u \cdot \nabla u = - \nabla p + \mu \Delta u \\
\mbox{div }u = 0\\
u \cdot \widehat{n} = 0 \mbox{ and } \omega = \mbox { curl }u = 0
\mbox{ at } \partial \Omega,
\end{array} \right.
\end{equation}

which is usually called the ``free boundary problem'' for the Navier
Stokes equations in the planar domain $\Omega$. The standard theory
for the Navier-Stokes equations with no-slip boundary condition can
be easily adapted to this problem, but even more, the vanishing
viscosity limit is quite well behaved in this case. Indeed, the
vorticity equation \eqref{vorteq} can be closed by adding the
Dirichlet boundary condition $\omega = 0$ at $\partial \Omega$ and
initial data. By multiplying the equation by $\omega^{p-1}$,
integrating, and performing the usual integration by parts, one
quickly arrives at the conclusion that the $L^p$-norms of vorticity
are bounded independently of $\nu$, and the compactness argument
outlined in the beginning of the section leads to an affirmative
answer to this variant of the $\partial$ Layer Problem. The
homogeneous Dirichlet condition is dissipative for the heat
equation; in particular, it means that the $L^p$ norms of vorticity
must decay in time.

There is a generalization of the free boundary condition called {\it
Navier friction condition} which is both physically and
mathematically interesting. The physically correct definition
involves the rate-of-strain matrix of the flow at the boundary but,
in two space dimensions, one can prove that this boundary condition
can be written as $\omega = (2\kappa - \alpha)u\cdot
\widehat{\tau}$, where $\kappa$ is the curvature of the boundary,
$\alpha$ is the friction coefficient and $\widehat{\tau}$ is the
counterclockwise unit tangent vector to $\partial \Omega$, if
$\Omega$ is a bounded domain with connected smooth boundary. We
refer the reader to \cite{CMR98} for a full discussion. This
condition, also called Navier ``slip'' condition, still allows for
the creation of vorticity at the boundary, but in a way that can be
controlled.

In \cite{CMR98}, T. Clopeau, A. Mikelic and R. Robert studied the
initial-boundary value problem for the Navier-Stokes equations in
2D, with Navier friction condition. They proved global
well-posedness for this problem for fixed viscosity, and they proved
the convergence of the vanishing viscosity limit to the (unique)
solution of the Euler equations if the initial vorticity is bounded.
An estimate on the rate of convergence was obtained by J. Kelliher
in \cite{kelliher06}.

In \cite{LNP05} the author, together with H. Nussenzveig Lopes and
G. Planas, extended the vanishing viscosity result by considering
flows with initial vorticity in $L^p$, $p>2$. The key point of this
result is the uniform $L^p$ estimate on vorticity, which is what we
will do next. Let $\Omega$ be a bounded, connected and
simply-connected smooth domain in $\real^2$ and $\alpha \in
C^{\infty}(\partial \Omega)$ . We consider $\omega^{\nu} =
\omega^{\nu}(x,t)$ to be the unique solution of the IBVP:

\begin{equation} \label{vorteq1}
\left\{ \begin{array}{l} \omega^{\nu}_t + u^{\nu} \cdot \nabla
\omega^{\nu} =
\nu\Delta \omega^{\nu}\\
u^{\nu} = K_{\Omega}[\omega^{\nu}]\\
\omega^{\nu}(x,t) = (2\kappa - \alpha)(u(x,t) \cdot \tau) \mbox{ for
} x \in \partial \Omega\\
\omega(x,0) = \omega_0
\end{array} \right.
\end{equation}

We assume that $\omega_0 \in H^1 \cap L^{\infty}$ is {\it
compatible}, i.e. that $\omega_0 = (2\kappa -
\alpha)(K_{\Omega}[\omega_0] \cdot \tau)$ in $\partial \Omega$.
Existence of a solution to problem \eqref{vorteq1} with
$\omega^{\nu} \in C([0,T);H^1(\Omega))$ (and $u^{\nu} \in
C([0,T);H^2(\Omega))$) was proved in \cite{CMR98}, and density of
compatible vorticities in $L^p$ with respect to the strong norm was
proved in \cite{LNP05}. Our objective is to prove the following

\begin{theorem} \label{LNP} (L., Nussenzveig Lopes and Planas, 2005)
For each $p>2$, there exists a constant $C>0$, independent of $\nu$,
such that
\[\|\omega^{\nu}(\cdot,t)\|_{L^p(\Omega)} \leq C (\|\omega_0
\|_{L^p(\Omega)} + \|u_0\|_{L^2(\Omega)}),\] with $u_0 =
K_{\Omega}[\omega_0]$.
\end{theorem}

\begin{proof}
In this proof we will use some results of the theory of parabolic
equations, which can be found in the book \cite{lieberman}.

The proof involves applying a maximum principle to two auxiliary
problems. First observe that $u\cdot \tau$ is $H^2$, and therefore
bounded. Set
\[\Lambda = \|(2\kappa - \alpha)u\cdot\tau\|_{L^{\infty}(\partial \Omega \times (0,T))}.\]
Consider the initial-boundary value problem for the Fokker-Planck
equation:

\begin{equation} \label{FP1}
\left \{ \begin{array}{ll}
 \widetilde{\omega}_t - \nu \Delta \widetilde{\omega} + u \cdot
\nabla \widetilde{\omega}=0 & \mbox{ in } \Omega \times (0,T),\\
\widetilde{\omega}(\cdot,0) = |\omega_0| & \mbox { in } \Omega,\\
\widetilde{\omega} = \Lambda & \mbox{ on }
\partial \Omega \times (0,T).
\end{array} \right.
\end{equation}

This problem has a unique weak solution $\widetilde{\omega} \in
L^2((0,T);H^1(\Omega))$. Then, $ \omega_1 = \omega -
\widetilde{\omega}$ is a weak solution for the following
initial-boundary value problem:

\begin{equation} \label{FP2} \left \{ \begin{array}{ll}
(\omega_1)_t - \nu \Delta \omega_1 + u \cdot
\nabla \omega_1 = 0 & \mbox{ in } \Omega \times (0,T),\\
\omega_1 (\cdot,0) = \omega_0 - |\omega_0| & \mbox { in } \Omega,\\
\omega_1 = ( 2 \kappa - \alpha ) u \cdot \tau - \Lambda & \mbox{ on
}
\partial \Omega \times (0,T).
\end{array} \right.
\end{equation}
The coefficients of the Fokker-Planck operator $\partial_t -\nu
\Delta + u \cdot \nabla$ are such that the maximum principle for
weak solutions of parabolic equations is valid. Therefore, as
$\omega_1 \leq 0$ on the parabolic boundary $\partial \Omega \times
(0,T)\cup \Omega \times \{t=0\}$, it follows that $\omega_1 \leq 0$
a.e. in $\Omega \times [0,T)$. Analogously, we prove that $\omega_2
= -\omega - \widetilde{\omega} $ is non-positive. We thus obtain

\begin{equation} \label{estimativa}
| \omega | \leq \widetilde{\omega} \mbox{ a.e. in } \Omega \times
[0,T).
\end{equation}

Moreover, since $\omega_0$ is bounded, the maximum principle may
also be applied to equation (\ref{FP1}) yielding that
$\widetilde{\omega} \in L^{\infty}((0,T)\times \Omega)$.

Next we obtain an estimate for $ \widetilde{\omega}. $ Let $
\widehat{\omega} = \widetilde{\omega} - \Lambda $. This is a
solution of the following problem:

\begin{equation} \label{widehateq}
 \left \{ \begin{array}{ll}
 \widehat{\omega}_t - \nu \Delta \widehat{\omega} + u \cdot
\nabla \widehat{\omega}=0 & \mbox{ in } \Omega \times (0,T),\\
\widehat{\omega}(\cdot,0) = |\omega_0| - \Lambda & \mbox { in } \Omega,\\
\widehat{\omega} = 0 & \mbox{ on }
\partial \Omega \times (0,T).
\end{array} \right.
\end{equation}

We formally multiply (\ref{widehateq}) by $ \widehat{\omega}
|\widehat{\omega}|^{p-2}$, where $ p > 2$, we integrate by parts and
use the incompressibility of the flow $u$ to obtain:
\begin{equation} \label{glu}
\frac{1}{p} \frac{d}{dt} \int_{\Omega} | \widehat{\omega} |^p +
(p-1) \nu \int_\Omega  \left||\nabla \widehat{\omega}
||\widehat{\omega}|^{(p-2)/2} \right| ^2 dx = 0.
\end{equation}
 Then,
\[
\| \widehat{\omega}(\cdot,t)\|_{L^p(\Omega)}  \leq  \|
\widehat{\omega}(\cdot,0)\|_{L^p(\Omega)}  \leq  \|\omega_0
\|_{L^p(\Omega)} + \Lambda |\Omega |^{1/p}. \] Therefore,
\begin{equation} \label{widetildest}
\| \widetilde{\omega} \|_{L^p(\Omega)} \leq \|
\widehat{\omega}\|_{L^p(\Omega)} + \Lambda |\Omega |^{1/p} \leq
\|\omega_0 \|_{L^p(\Omega)} + 2 \Lambda |\Omega |^{1/p}.
\end{equation}

This formal calculation can be made rigorous by using the weak
formulation of (\ref{widehateq}) given in \cite{lieberman}. We begin
by observing that $\widehat{\omega}_t \in L^2((0,T);H^{-1}(\Omega))$
and $\widehat{\omega} \in L^2((0,T);H^1_0(\Omega)) \cap
L^{\infty}((0,T)\times \Omega )$. This implies that
$\widehat{\omega} |\widehat{\omega}|^{p-2} \in
L^2((0,T);H^1_0(\Omega))$. Therefore we can multiply
(\ref{widehateq}) by $\widehat{\omega} |\widehat{\omega}|^{p-2}$ if
we understand the product with $\widehat{\omega}_t$ and with $\Delta
\widehat{\omega}$ as duality pairings. Finally, in order to justify
\eqref{glu} one still needs to approximate $\widehat{\omega}$ by
suitable smooth functions and pass to the limit in each term of the
weak formulation so as to obtain
\[\frac{1}{p}\frac{d}{dt}\int_{\Omega}|\widehat{\omega}|^p = \langle \widehat{\omega}_t ,
\widehat{\omega} |\widehat{\omega}|^{p-2} \rangle \]

\[=\nu \langle \Delta \widehat{\omega},\widehat{\omega} |\widehat{\omega}|^{p-2} \rangle = -(p-1) \nu
\int_\Omega  ||\nabla \widehat{\omega}
||\widehat{\omega}|^{(p-2)/2} | ^2 dx. \]

This can be easily accomplished using mollification in time together
with the Dirichlet heat semigroup for $\Omega$, thus generating a
family of smooth functions $\widehat{\omega}_{\varepsilon}$ such
that $\partial_t \widehat{\omega}_{\varepsilon} \to
\widehat{\omega}_t$ strongly in $ L^2((0,T);H^{-1}(\Omega))$, while
$\widehat{\omega}_{\varepsilon}
|\widehat{\omega}_{\varepsilon}|^{p-2} \rightharpoonup
\widehat{\omega} |\widehat{\omega}|^{p-2}$ weakly in $
L^2((0,T);H^1_0(\Omega))$ and $\widehat{\omega}_{\varepsilon}$ is
uniformly bounded in $\Omega \times (0,T)$.

Given (\ref{widetildest}) we now turn to the estimate of $\Lambda$.
Using Sobolev imbedding and interpolating between $W^{1,p}$ and
$L^2$, we find:
\[\| u^{\nu} (\cdot,t) \cdot \tau \|_{L^{\infty}(\partial \Omega)} \leq C \|u^{\nu} (\cdot,t) \|_{C(\bar{\Omega})}
\leq C \|u^{\nu} (\cdot,t)\|_{L^2(\Omega)}^\theta \| u^{\nu}
(\cdot,t) \|_{W^{1,p}(\Omega)}^{1- \theta} \]
\[\leq C \|u^{\nu} (\cdot,t) \|_{L^2(\Omega)}^\theta \| \omega^{\nu} (\cdot,t)\|_{L^p(\Omega)}^{1- \theta},\]
where $\theta = (p-2)/(2p-2)$.

Let $\varepsilon$ be an arbitrary positive number. We now use
Young's inequality together with the fact that $\kappa$ and $\alpha$
are bounded to conclude that:
\begin{equation} \label{estimativa2}
\Lambda \leq C_{\varepsilon}  \|u^{\nu}
\|_{L^{\infty}((0,T);L^2(\Omega)} +  \varepsilon \| \omega^{\nu}
\|_{L^{\infty}((0,T);L^p(\Omega))}
\end{equation} for some $C_{\varepsilon}>0$. Taking $\varepsilon$
small enough, from (\ref{estimativa})-(\ref{estimativa2}) we obtain:
\begin{equation} \label{estimatevort}
\| \omega^{\nu} \|_{L^\infty(0,T;L^p(\Omega))}  \leq C( \|\omega_0
\|_{L^p(\Omega)}  + \|u^{\nu} \|_{L^\infty(0,T;L^2(\Omega))})
\end{equation} for any $p > 2$, where $ C =
C(p,\Omega,\|\kappa\|_{L^{\infty}(\partial
\Omega)},\|\alpha\|_{L^\infty(\partial \Omega)} ).$ Finally, a
standard energy estimate, yields $\|u \|_{L^\infty(0,T;L^2(\Omega))}
\leq \|u_0\|_{L^2(\Omega)}$, thereby concluding the proof.

\end{proof}

Let us first observe that the "standard energy estimate" used as
last step of the proof above is far from obvious. The difficulty
lies in the integration by parts of the viscous term in the energy
estimate. This can actually be done by using the correct formulation
of the Navier friction condition, namely $2(DU)_s\widehat{n} \cdot
\widehat{\tau} + \alpha u \cdot \widehat{\tau} = 0$, where $(DU)_s$
is the symmetric part of the Jacobian matrix $DU$. We leave the
derivation of the energy estimate in this case as an exercise to the
reader.

The estimate in Theorem \ref{LNP} can be extended to any vorticity
in $L^p$, $p>2$ by density. It is not clear whether the restriction
$p>2$ is a consequence of the technique used above, or it is an
essential restriction in the estimate above. This is an interesting
open problem. One physical justification for the Navier friction
boundary condition is that is approximates the interaction of
incompressible flow with a rough boundary. This contention was
rigorously justified by Jager and Mikelic, who derived the Navier
friction condition as the effective behavior associated with the
homogenization of an oscillatory boundary, see \cite{jm01}.

There is no physical interpretation, at present, for the vanishing
viscosity limit associated with Navier friction condition. It is,
however, a very natural question from the mathematical point of
view, since it interpolates naturally between the trivial case of
the free boundary condition and the Boundary Layer Problem. The
Navier friction condition still gives rise to a boundary layer,
albeit in a less singular way than Prandtl's original theory. This
is the subject of the next Section.

\section{Boundary layer for the Navier friction condition}

In this section we will explore a rigorous procedure to obtain a
boundary layer equation, based on ideas from geometric optics,
applied to the Navier friction condition. What we will present here
is based on joint work of the author with D. Iftimie, H. Nussenzveig
Lopes and F. Sueur.

We begin with the $\nu$-Navier-Stokes system

\begin{equation} \label{nuNS}
\left\{ \begin{array}{l}
\partial_t u^{\nu} + u^{\nu} \cdot \nabla u^{\nu} = -\nabla p^{\nu} + \nu \Delta u^{\nu}\\ \\
\mbox{div } u^{\nu} = 0\\ \\
u^{\nu}_2(t,x_1,0) = 0\\ \\
\omega^{\nu}(t,x_1,0) = \gamma u^{\nu}_1(t,x_1,0)\\ \\
u^{\nu}(0,x) = u_0(x). \end{array} \right.
\end{equation}

We denote by $u^E$, $\omega^E$, $p^E$ the solution of the
initial-boundary value problem for the Euler equations in the
half-plane, with initial data $u_0(x)$ and the standard
no-penetration boundary condition.

We set

\[u^{\nu} = u^E + \sqrt{\nu} \left( \begin{array}{l} v(t,x_1,x_2/\sqrt{\nu})\\ 0 \end{array} \right)
+ \nu u^{\nu}_R. \]

Additionally, we write

\[\omega^{\nu} = \omega^E + W(t,x_1,x_2/\sqrt{\nu}) + \sqrt{\nu} \phi^{\nu},\]
where $W = W(t,x_1,y) = (\partial_y v)(t,x_1,y)$ and $\phi^{\nu} =
\sqrt{\nu} \mbox{ curl }u^{\nu}_R$.

We are going to write an equation for $v$, derive an equation for
$W$ and then use the equations for $W$, $\omega^{\nu}$ and
$\omega^E$ to deduce an equation for $\phi^{\nu}$. Our final
objective is to obtain an $L^p$ estimate for $\phi^{\nu}$ which is
independent of $\nu$.

We first write  an evolution equation with Neumann boundary
condition for the velocity profile of the boundary layer $v$:

\begin{equation} \label{veq}
\left\{ \begin{array}{l}
\partial_t v + u_1^E(t,x_1,0) \partial_{x_1}v + \partial_{x_1}u^E_1(t,x_1,0) v
+ y (\partial_{x_2} u_2^E)(t,x_1,0) \partial_y v = \partial_y^2 v\\ \\
v(0,x_1,y) = 0\\ \\
\partial_y v(t,x_1,0) = \gamma u_1^E(t,x_1,0) - \omega^E(t,x_1,0). \end{array}
\right.
\end{equation}

Now we differentiate equation \eqref{veq} with respect to $y$ and
use the div-free condition on $u^E$ to obtain an equation for the
vorticity profile of the boundary layer $W$. We find:

\begin{equation} \label{Weq}
\left\{ \begin{array}{l}
\partial_t W + u_1^E(t,x_1,0) \partial_{x_1}W + y (\partial_{x_2} u_2^E)(t,x_1,0) \partial_y W
= \partial_y^2 W\\ \\
W(0,x_1,y) = 0\\ \\
W(t,x_1,0) = \gamma u_1^E(t,x_1,0) - \omega^E(t,x_1,0). \end{array}
\right.
\end{equation}

The equations for $v$ and $W$ are linear, with smooth coefficients
and independent of $\nu$. For the sake of these lectures, we assume
well-posedness and exponential decay of the solutions $v$ and $W$ as
$y \to \infty$. To be more precise, we assume existence of a unique
smooth solution $v=v(t,x_1,y)$ of \eqref{veq} and consequently of $W
= W(t,x_1,y)$, solution of \eqref{Weq}, such that:
\begin{enumerate}
\item $W \in C^{\infty}_b \cap L^p$, for some $p>2$,
\item $v \in C^{\infty}_b \cap L^p$ and $\partial_{x_1} v \in L^p$,
\item $\partial_{x_1} W$ and $\partial_{x_1}^2 W$ both in $L^p$,
\item $y \partial_y W$  bounded and $y^2 \partial_y W \in L^p$,
\item $y \partial_{x_1} W \in L^p$.
\end{enumerate}

Next we write the equation for $\phi^{\nu}$. First set $z =
x_2/\sqrt{\nu}$. We introduce

\[G^{\nu} = [u^E_1(t,x_1,x_2) - u^E_1(t,x_1,0)] (\partial_{x_1} W)(t,x_1,z) \]
\[ + [u^E_2(t,x_1,x_2) - x_2 (\partial_{x_2}u^E_2)(t,x_1,0)] \frac{1}{\sqrt{\nu}} (\partial_y W)(t,x_1,z), \]

and

\[F^{\nu} = \sqrt{\nu} \Delta \omega^E + \sqrt{\nu} (\partial_{x_1}^2 W)(t,x_1,z)
- \left[ \left( \begin{array}{l} v(t,x_1,z)\\0 \end{array} \right) +
\sqrt{\nu} u^{\nu}_R \right] \cdot \nabla \omega^E\]
\[- (v(t,x_1,z) + \sqrt{\nu} (u^{\nu}_R)_1) (\partial_{x_1}W)(t,x_1,z)  -
(u^{\nu}_R)_2 (\partial_y W)(t,x_1,z) - \frac{G^{\nu}}{\sqrt{\nu}}.
\]

With this notation, the equation for $\phi^{\nu}$ becomes

\begin{equation} \label{phinueq}
\left\{ \begin{array}{l}
\partial_t \phi^{\nu} + u^{\nu} \cdot \nabla \phi^{\nu} - \nu \Delta \phi^{\nu} = F^{\nu}\\ \\
\phi^{\nu}(0,x) = 0 \\ \\
\phi^{\nu}(t,x_1,0) = \gamma(v+\sqrt{\nu}(u^{\nu}_R)_1)(t,x_1,0).
\end{array} \right.
\end{equation}

We break problem \eqref{phinueq} into two problems, one with
homogeneous boundary condition and nonzero forcing and another with
nonhomogeneous boundary data and zero forcing. We write $\phi^{\nu}
\equiv \phi^{\nu}_1 + \phi^{\nu}_2$, with

\begin{equation} \label{phinu1eq}
\left\{ \begin{array}{l}
\partial_t \phi^{\nu}_1 + u^{\nu} \cdot \nabla \phi^{\nu}_1 - \nu \Delta \phi^{\nu}_1 = 0\\ \\
\phi^{\nu}_1(0,x) = 0 \\ \\
\phi^{\nu}_1(t,x_1,0) = \gamma(v+\sqrt{\nu}(u^{\nu}_R)_1)(t,x_1,0),
\end{array} \right.
\end{equation}

and

\begin{equation} \label{phinu2eq}
\left\{ \begin{array}{l}
\partial_t \phi^{\nu}_2 + u^{\nu} \cdot \nabla \phi^{\nu}_2 - \nu \Delta \phi^{\nu}_2 = F^{\nu}\\ \\
\phi^{\nu}_2(0,x) = 0 \\ \\
\phi^{\nu}_2(t,x_1,0) = 0. \end{array} \right.
\end{equation}

The basic idea is to treat \eqref{phinu2eq} using the energy method
and to treat \eqref{phinu1eq} adapting the comparison principle
argument used in \cite{LNP05}. In addition, it is necessary to
obtain a priori estimates for the elliptic system which relates
$\phi^{\nu}$ with $u^{\nu}_R$. We write this system as:

\begin{equation} \label{ellip}
\left\{ \begin{array}{l} \mbox{div } \sqrt{\nu} u^{\nu}_R =
\partial_{x_1}v(t,x_1,x_2/\sqrt{\nu}) \\
\mbox{curl }\sqrt{\nu} u^{\nu}_R = \phi^{\nu} \\
(u^{\nu}_R)_2(t,x_1,0) = 0.
\end{array} \right.
\end{equation}

What we end up  with is the following estimates
\[\|\phi^{\nu}_1\|_{L^p} \leq C(1+\|\phi^{\nu}_2\|_{L^p}),\]
and
\[\frac{d}{dt} \|\phi^{\nu}_2\|_{L^p} \leq C (\|\phi^{\nu}\|_{L^p} + 1) \leq
C( 1 + \|\phi^{\nu}_1\|_{L^p} + \|\phi^{\nu}_2\|_{L^p}),\] the
proofs are the technical core of this work, but for time limitations
we omit them.

Just to conclude this story, the estimates above lead to
\[\frac{d}{dt} \|\phi^{\nu}_2\|_{L^p} \leq C(1 + \|\phi^{\nu}_2\|_{L^p}). \]
It follows from Gronwall's inequality that $\|\phi^{\nu}_2\|_{L^p}$
is bounded independently of $\nu$, for a fixed time interval
$[0,T]$. Then, using the first estimate again concludes that there
exists a constant $C$ independent of $\nu$ such that
\[\|\phi^{\nu}\|_{L^p} \leq C. \]

{\bf Remark:} We observe that in \cite{is08}, Iftimie and Sueur
proved that $\phi^{\nu}$ is bounded in $L^2((0,T)\times\Omega)$,
independently of $\nu$. We proved that there is a bound on
$L^{\infty}((0,T);L^p)$, for $2<p<\infty$. We do not expect to be
able to extend this argument to $p = \infty$. Note also that the
analysis in \cite{is08} is valid for general smooth domains in two
or three space dimensions, whereas the analysis outlined above only
works in two dimensions.

\section{The rotating cylinder}

Our purpose in this section is to provide a fairly explicit
illustration of the boundary layer phenomena. The main point of this
illustration is to show, that in a favorable scenario, the small
viscosity limit has vorticity behaving like a vortex sheet. This
explains, in part, the difficulty involved in the Boundary Layer
Problem, since vortex sheet regularity is critical for the weak
continuity of the nonlinearity present in Euler and Navier-Stokes.
This material is based on joint work of the author with A.
Mazzucato, H. Nussenzveig Lopes and M. Taylor, and it is contained
in two articles, \cite{LMN08,LMNT07}.

The physical situation we wish to consider is that of an infinite
circular cylinder filled with fluid. The boundary of the cylinder is
a material shell which is assumed to rotate about its center of
symmetry with angular velocity $\alpha = \alpha(t)$ (positive
rotation is counterclockwise). We restrict our attention to planar
and circularly symmetric motions, so that the fluid velocity $u$ is
given by
\[u = u(x_1,x_2,x_3,t) = v(r,t)\frac{(-x_2,x_1,0)}{r} = V(x_1,x_2,t)\frac{(-x_2,x_1,0)}{r},\]
with $r=\sqrt{x_1^2+x_2^2}$. This symmetry assumption is consistent
with the Navier-Stokes equations as long as the initial velocity is
planar and circularly symmetric, i.e., weak solutions satisfying the
symmetry assumptions exist globally in time for symmetric initial
data.  It should be remarked that the assumption that such flows
stay planar for all Reynolds number is unphysical - 3D turbulence is
expected if the Reynolds number is large enough. The limit behavior
we wish to present is motivated by mathematical considerations.

Let us first observe that, if velocity is circularly symmetric, and,
in particular, tangent to concentric circles around the origin,
vorticity is also circularly symmetric, and therefore constant in
the same circles. This implies that $\nabla \omega$ is perpendicular
to these circles, and therefore $u$ and $\nabla \omega$ are
orthogonal everywhere. Therefore, the nonlinearity $u \cdot \nabla
\omega$ is the vorticity equation vanishes identically, and the
vorticity equation becomes the heat equation. This basic fact is
what it makes possible to prove the following result.

\begin{theorem} \label{thm:maintheorem}
Let $u^{\nu}$ be the solution of the 2D Navier-Stokes  equations in
the unit disk, with no slip boundary data with respect to boundary
rotation with prescribed angular velocity $\alpha \in BV(\real)$.
Assume that the initial velocity $u_0\in L^2(\disk)$ has circular
symmetry, i.e. $u_0(x) = v_0(|x|) x^{\perp}$. Then, $u_0$ is a
steady solution of the 2D Euler equation and $u^{\nu}$ converges
strongly in $L^{\infty}([0,T], L^2(\disk))$ to $u_0$ as $\nu \to 0$.
\end{theorem}

The proof of this result is based on a semigroup treatment of a
symmetry-reduced problem. The important issue is that the symmetry
assumption not only prevents boundary layer separation, but it also
eliminates the nonlinearity of the problem. Thus, the
symmetry-reduced problem is a linear, singular coefficient
perturbation of the heat equation, and its treatment is classical.

Let $D \equiv \{|x| < 1\} \subset \real^2$. We begin by considering
the following $2\times 2$ system of PDE, consisting of a pair of
uncoupled heat equations:

\begin{equation} \label{circsymNS}
\left\{\begin{array}{ll}
\partial_t u^{\nu} = \nu \Delta u^{\nu}, & \mbox{ in } (0,\infty)\times D;\\
u^{\nu}(x,0) = u_0(x), & \mbox{ in } D;\\
u^{\nu}(x,t) = \frac{\alpha(t)}{2\pi} x^{\perp}, & \mbox{ on }
|x|=1.
\end{array}\right.
\end{equation}

Let us suppose that $u_0$ is {\it covariant under rotations}, i.e.,
$u_0(R_{\theta}x)=R_{\theta}u_0(x)$ for any rotation $R_{\theta}$.
Then \eqref{circsymNS} is the symmetry-reduced 2D Navier-Stokes
system of equations with initial velocity $u_0$.

Let us assume, to begin with, that $u_0\in C^{\infty}(\overline{D})$
is covariant under rotations and that $u_0$ and $\alpha \in
C^{\infty}(\real)$ satisfy the following compatibility conditions:

\begin{equation} \label{compcond}
\begin{array}{l}
u_0(x)= \alpha(0)\frac{x^{\perp}}{2\pi} \mbox{ on } \partial D; \\
\\
\Delta u_0 = \alpha^{\prime} (0)\frac{x^{\perp}}{2\pi\nu} \mbox{ on
} \partial D.
\end{array}
\end{equation}

Then it follows that there exists a unique weak solution $u^{\nu} =
u^{\nu}(x,t) \in L^2(\real;H^4(D)) \cap H^1(\real;H^2(D)) \cap
H^2(\real;L^2(D))$ of \eqref{circsymNS}  In particular we have
$u^{\nu} \in C^0(\real;C^1(\overline{D}))$. With this we can now
establish the following lemma.

\begin{lemma} \label{Neumannreg}
 Let $\omega^{\nu} = \nabla^{\perp} \cdot u^{\nu} \equiv \partial_{x_1} u^{\nu}_2 - \partial_{x_2}u^{\nu}_1$.  Then
$\omega^{\nu}$ is a classical solution of
\begin{equation} \label{circsymvortNS}
\left\{\begin{array}{ll}
\partial_t \omega^{\nu} = \nu \Delta \omega^{\nu}, & \mbox{ in } (0,\infty)\times D;\\
\frac{\partial \omega^{\nu}}{\partial \widehat{\mathbf{n}}}(x,t) = \frac{\alpha^{\prime}(t)}{2\pi \nu}, & \mbox{ on }    \{|x|=1\} \times [0,\infty)\\
\omega^{\nu}(x,0) = \nabla^{\perp} \cdot u_0(x), & \mbox{ in } D
\times \{t=0\}.
\end{array}\right.
\end{equation}
Additionally, $\int_D \omega^{\nu} (x,t)\,dx = \alpha(t)$ at every
$t > 0$.
\end{lemma}

\begin{proof}
Start by noting that, if $u_0$ is covariant under rotations about
the origin, then $\omega_0$ is circularly symmetric and, hence, so
is $\omega^{\nu}(\cdot,t)$; this is due to the rotational invariance
of the heat equation.

The partial differential equation in \eqref{circsymvortNS} can be
trivially deduced from \eqref{circsymNS}, along with the initial
condition, so that all that remains is to verify that the boundary
condition is correct.

We integrate \eqref{circsymvortNS} in space and use the divergence
theorem to obtain:
\[\frac{d}{dt}  \int_{D} \omega^{\nu}(x,t)\,dx
      = \nu \int_{\partial D} \nabla \omega^{\nu}\cdot x \,dS\]
      \[= 2\pi \nu \frac{\partial \omega^{\nu}}{\partial \widehat{\mathbf{n}}}\lfloor{\partial D},\]
where we used the circular symmetry of $\omega^{\nu}$ in the last
step.

Next observe that $\omega^{\nu} = -\mbox{ div }(u^{\nu})^{\perp}$,
so that
\[\int_D \omega^{\nu} (x,t)\, dx = -\int_{\partial D} (u^{\nu})^{\perp} \cdot \mathbf{n} \, dS = \alpha (t),\]
as $u^{\nu}$ is a solution of \eqref{circsymNS}.

These two facts yield the desired conclusions.

\end{proof}

Let us note in passing that the identity $\int_D \omega^{\nu}
(x,t)\,dx = \alpha(t)$ is valid (a.e. in time) under the much weaker
assumption $u_0 \in L^2(D)$ and $\alpha \in BV(\real)$, without the
compatibility conditions \eqref{compcond}, since then the solution
$u^{\nu} \in L^2((0,T);H^1(D))$.

Next we state and prove our main theorem in which we examine the
inviscid limit of $\omega^{\nu}$.

\begin{theorem}

Let $\alpha \in  BV(\real)$ and assume that $\alpha$ is compactly
supported in $(0,\infty)$. Let $u_0 \in L^2(D;\real^2)$ be covariant
under rotations and assume that $\omega_0 = \nabla^{\perp} \cdot u_0
\in L^1(D)$. Then we have:
\begin{enumerate}
\item There exists a constant $C>0$ such that
\[\int_D |\omega^{\nu} (x,t)|\,dx \leq C(\|\omega_0\|_{L^1} + \|\alpha\|_{BV}),\]
for almost all time.
\item For almost every $t \in \real$,
\[\alpha (t) = \int_D \omega^{\nu}(x,t) \, dx.\]
\item We also have, for any $0<a<1$,
\[\|\omega^{\nu} - \omega_0\|_{L^1(\{|x|\leq a\})} \to 0 \mbox{ as } \nu \to 0,\]
a.e. in time. If $\omega_0 \in C^0(D)$ then the convergence is
uniform in time.
\end{enumerate}
\end{theorem}

Given $\alpha \in BV(\real)$, we call the modified $BV$ function
$\widetilde{\alpha}$ the left-continuous ``correction" of $\alpha$.

\begin{proof}
We start by choosing appropriate regularizations of $\alpha$ and of
$u_0$. Let $\alpha_n$ be the regularization obtained in Lemma 4.1 of
\cite{LMN08}, satisfying the following properties:
\begin{enumerate}
\item $\alpha_n \in C^{\infty}(\real)$;
\item $\alpha_n \to \widetilde{\alpha}$ pointwise and strongly in $L^p_{\loc}(\real)$, for any $1\leq p < \infty$;
\item $\alpha_n^{\prime} \rightharpoonup \alpha^{\prime}$ weak-$\ast$ in $\mathcal{BM}(\real)$.
\end{enumerate}
Above, the derivatives are understood in the sense of distributions.
Let $u_{0,n}$ be a sequence of vector fields in $C^{\infty}_c(D)$,
commuting with rotations, such that $u_{0,n} \to u_0$ strongly in
$L^2(D)$. Since $\alpha$ is compactly supported in $(0,\infty)$ it
follows immediately that the compatibility conditions
\eqref{compcond} are satisfied for $u_{0,n}$ and $\alpha_n$. Hence,
if $u^{\nu}_n$ denotes the solution of the symmetry-reduced
Navier-Stokes equations \eqref{circsymNS} with data $u_{0,n}$ and
$\alpha_n$ then $\omega^{\nu}_n= \nabla^{\perp} \cdot u^{\nu}_n$ is
a solution of \eqref{circsymvortNS} such that
\begin{equation} \label{hah}
\alpha_n(t) = \int_D \omega^{\nu}_n(\cdot,t),
\end{equation}
by Lemma \ref{Neumannreg}.

Let us prove the estimate in item (1) for $\omega^{\nu}_n$.

We will derive an energy-type estimate in $L^1$ after appropriately
regularizing $|\omega^{\nu}_n |(x,t)= sgn(\omega^{\nu}_n
(x,t))\omega^{\nu}_n (x,t)$. We first mollify $y \mapsto |y|$. For
each $\vare
> 0$ there exists a $C^{1}$, convex function $\phi_{\vare}$, such
that $|\phi'_{\vare}|$ is bounded, uniformly in $\vare$, and
$\phi_{\vare} \to |\cdot|$ as $\vare \to 0$.

Next, we multiply the PDE in \eqref{circsymvortNS} by
$\phi'_{\vare}(\omega^{\nu}_n)$ and  integrate over the disk to
obtain:
\[
    \frac{d}{dt}  \int _{|x|\leq 1}\phi_{\vare}(\omega^{\nu}_n) \,dx =
     \nu \int_{x|\leq 1} \left( \Delta[\phi_{\vare}(\omega^{\nu}_n)]
     - \phi''_{\vare}(\omega^{\nu}_n)|\nabla \omega^{\nu}_n |^2\right)\,dx.
\]
Since $\phi_{\vare}\in C^1$ and convex, $\phi''_{\vare}
(\omega^{\nu}_n) \geq 0$ so that, by the divergence theorem,
\[
  \begin{aligned}
   \frac{d}{dt}  \int _{|x|\leq 1}\phi_{\vare}(\omega^{\nu}_n) \,dx
   &\leq \nu \int_{|x|\leq 1} \Delta[\phi_{\vare}(\omega^{\nu}_n)]\,dx =
   \nu \int_{|x|=1} \nabla [\phi_{\vare}(\omega^{\nu}_n)]  \cdot x \,
    dS(x)  \\
   &= \frac{\alpha_n'(t)}{2\pi} \int_{|x|=1} \phi'_{\vare}(\omega^{\nu}_n)
   \,dS(x) \leq C\,|\alpha_n'(t)|,
  \end{aligned}
\]
where we used the boundary condition for $\omega^{\nu}_n$ in
\eqref{circsymvortNS}. Then, integrating in time gives:
\[
     \int_{|x|\leq 1} \phi_{\vare}(\omega^{\nu}_n)\,dx  \leq
      \int_0^t |\alpha_n'(\tau)|\, d\tau + \int_{|x|\leq 1}
      \phi_{\vare} (\omega_{0,n})\,dx.
\]
Observe that we can choose $\phi_{\vare}$ in such a way that
$\phi_{\vare}(y)\to |y| $ monotonically. Consequently, by the
Monotone Convergence Theorem, passing to the limit $\vare \to 0$ in
the
 expression above gives, for every $0<t < \infty$:
\begin{equation} \label{approximest}
        \int_{|x|\leq 1} |\omega^{\nu}_n(x,t)|\,dx \leq
      \int_0^t |\alpha_n'(\tau)|\, d\tau + \int_{|x|\leq 1}  |\omega_{0,n}(x)|\,dx \leq
\|\alpha\|_{BV} + \|\omega_0\|_{L^1},
\end{equation}
from which we obtain item (1) {\bf for } $\mathbf{\omega^{\nu}_n}$.

We have shown that, for each $\nu>0$,
$\{\omega^{\nu}_n\}_{n=1}^{\infty}$ is uniformly bounded in
$L^{\infty}(\real;L^1(D))$. Thus, passing to subsequences as needed,
we find that $\{\omega^{\nu}_n\}_{n=1}^{\infty}$ is weak-$\ast$
$L^{\infty}(\real;\mathcal{BM}(D))$ convergent as $n \to \infty$.

On the other hand, the properties of $\alpha_n$, together with the
proofs of Proposition 4.3 and Proposition 5.1 in \cite{LMN08} imply
that $u^{\nu}_n \rightharpoonup u^{\nu}$ weakly in $L^2((0,T)\times
D)$, where $u^{\nu}$ is the solution of \eqref{circsymNS}, with data
$u_0$ and $\alpha$. It follows immediately that
\[\omega^{\nu}_n \; \rightharpoonup \; \omega^{\nu} = \nabla^{\perp}\cdot u^{\nu}, \]
weakly in $L^2((0,T); H^{-1}(D))$. Uniqueness of limits implies that
the convergence of $\omega^{\nu}_n$ to $\omega^{\nu}$ can be
improved to weak-$\ast$ $L^{\infty}(\real;\mathcal{BM}(D))$. We use
again an estimate obtained in the proof of Proposition 5.1 of
\cite{LMN08}, see (5.3) there, to find that $\{\omega^{\nu}_n\}_n$
is bounded uniformly in $L^2((0,T) \times D)$. Hence we obtain
\[\omega^{\nu}_n \rightharpoonup \omega^{\nu} \mbox{ weakly in } L^2((0,T) \times D). \]

Therefore, since $L^2((0,T) \times D) \subset L^2((0,T);L^1(D))$, it
follows that $\omega^{\nu} (\cdot,t) \in L^1(D)$ a.e. $t>0$. Weak
lower-semicontinuity of the $L^1$-norm implies that item (1) follows
from \eqref{approximest}.

Since $\alpha_n \to \alpha$ strongly in $L^2_{\loc}(\real)$ it
follows that, for any $\varphi = \varphi(x) \in L^{\infty}(D)$ and
$\psi = \psi(t) \in L^2((0,T))$ we have
\[\int_0^T\int_D \left( \frac{\alpha_n(t)}{\pi} - \omega^{\nu}_n(x,t)\right)\varphi(x)\psi(t) \,dxdt \to
\int_0^T\int_D \left(\frac{\alpha(t)}{\pi} -
\omega^{\nu}(x,t)\right)\varphi(x)\psi(t) \,dxdt.\] Take $\varphi(x)
\equiv 1$. Then, using \eqref{hah}, we see that the left-hand-side
above vanishes identically, so that
\[\int_0^T\int_D \left(\frac{\alpha(t)}{\pi} - \omega^{\nu}(x,t)\right)\psi(t) \,dxdt = 0,\]
i.e.
\[\int_0^T \left(\alpha(t) - \int_D \omega^{\nu}(x,t)\,dx\right)\psi(t)\,dt = 0\]
for any $\psi \in L^2((0,T))$. Take $\psi=\psi(t)=\alpha(t) - \int_D
\omega^{\nu}(x,t)\,dx$; this gives item (2).

Note that $\omega^{\nu}$ is a solution of the heat equation with
viscosity $\nu$ in $D$, with initial data $\omega_0$.

To prove (3), let $\widetilde{\Omega}$ be an open, compactly
contained subset of $D$ and let $\phi \in C^{\infty}_c(D)$ be such
that $\phi \equiv 1$ in a neighborhood of $\widetilde{\Omega}$. We
consider $v^{\nu} \equiv \phi \omega^{\nu}$, extended to the full
plane. Then $v^{\nu}$ satisfies:
\[\partial_t v^{\nu} = \nu \Delta v^{\nu} + F^{\nu},\]
with
\[F^{\nu} = - \nu ( 2 \nabla \phi \cdot \nabla \omega^{\nu} + \omega^{\nu} \Delta \phi).\]

We apply Duhamel's formula in the whole plane to obtain:

\begin{equation} \label{gabby}
v^{\nu} = e^{\nu t \Delta}(\phi \omega_0) + \int_0^t e^{\nu (t-s)
\Delta} F^{\nu}(s) \, ds.
\end{equation}

We wish to estimate $v^{\nu}=v^{\nu}(x,t)$ for $x \in
\widetilde{\Omega}$.

Clearly, the first term converges to $\omega_0$ in
$\widetilde{\Omega}$, as $\nu \to 0$, in whichever space $\omega_0$
lies in. What we are left to prove is that the other term in
\eqref{gabby} converges to zero. This follows from the fact that the
heat kernel is sharply localized when $\nu \to 0$, together with the
fact that the support of $F^{\nu}$ is bounded away from
$\widetilde{\Omega}$, uniformly in $\nu$, as we are only interested
in $x \in \widetilde{\Omega}$.

\end{proof}

From this result we derive three conclusions. The first is that, if
$\alpha \in BV(\real)$ then $\omega^{\nu} \rightharpoonup W$
weak-$\ast$ in $L^{\infty}_{\loc}((0,\infty);BM(\overline{D}))$ and
hence $W = \omega_0 + \mu$, where $\mu$ is a measure supported on
the boundary of the disk $D$. This implies that, in considering the
classical open problem of the inviscid limit for the Navier-Stokes
equations in domains with boundary, at the very least one has to
deal analytically with regularity at the level of vortex sheets,
without a priori sign conditions.  The second observation is that
the vorticity generated by a boundary layer appears to be
proportional to the {\it acceleration} of the boundary with respect
to the adjoining flow.  Third, we have proved that, if $\alpha$ is
not constant, then the vorticity $\omega^{\nu}$ does not converge in
$L^1$ to the vorticity of the limit flow $\omega_0$. This shows that
the $L^2$ convergence of velocity fields obtained in \cite{LMN08}
cannot be improved to convergence in derivatives. Finally, we would
like to mention that, in \cite{LMNT07}, we weaken the hypothesis
that $\alpha \in BV$ to $\alpha \in L^p$, for $p>1$, and we adjust
the convergence obtained accordingly. However, the vorticity
accounting we have shown here only works for $\alpha \in BV$.

\section{Conclusion}

 As a conclusion for these notes, it would be interesting
 to highlight a few main points of the our discussion on the
 mathematical theory of boundary layers.

 \begin{enumerate}

 \item We remark on the extent of the current lack of
 understanding of the vanishing viscosity limit in the presence of
 boundaries. At this stage, the Boundary Layer Problem formulated
 as formulated in Section 4 is wide open.

 \item The Prandtl equation and its refinements take a magnifying
 glass approach to the difficulty involved in boundary layers. This
 approach is very interesting and useful for specific problems, but the
 difficulties in the well-posedness of Prandtl's equation do not make it
 very promising as an avenue for the Boundary Layer Problem.

 \item The magnifying glass approach works well for the case of
 Navier friction condition mainly because the boundary layer
 equation in this case is linear and well-posed.

 \item At best, the Boundary Layer Problem involves treating a
 sequence of approximate solutions to the Euler equations with
 vortex-sheet level regularity, which is the critical case for
 the weak continuity of the nonlinearity.

 \item The strength of the boundary layer as a vortex sheet, and
 therefore some of the difficulty in the Boundary Layer Problem
 appears to be associated with the acceleration of fluid past the
 boundary.

\end{enumerate}

Finally let us put together a few problems connected with what was
presented here. First, establish rigorously the ill-posedness of
Prandtl's equation without monotonicity. Second, we mention
extending the $L^p$ control on the vorticity in the case of Navier
friction condition to $p \leq 2$. Third, in \cite{LNX01}, the
author, together with H. Nussenzveig Lopes and Zhouping Xin
introduced the notion of {\it boundary coupled weak solution} of the
incompressible Euler equations, and we proved the existence of such
a weak solution in the case of the half-plane. It would be
interesting to know if the vanishing viscosity limit, for example,
for Navier friction condition in the half-plane, gives rise to such
weak solutions. Also, if vorticity is $L^p$, $p\geq 2$ the theory of
renormalized solutions of DiPerna-Lions implies that weak solutions
of the Euler equations conserve $p$-norms exactly, for flows in
domains without boundary. Does this remain true in the case of
vanishing viscosity limits in domains with boundary? Would it be
possible to find an example, in the spirit of the circularly
symmetric flows in Section 7, but for which there is boundary layer
separation? One possibility is to look for the solution of the
Navier-Stokes system on the disk, with initial vorticity given by an
odd eigenfunction of the Dirichlet Laplacian. Finally, is boundary
layer separation possible for flows with Navier friction condition?
A related problem is to study the vanishing viscosity limit with
Navier friction condition in nonsmooth domains. Another pair of
problems is to extend the boundary layer analysis done in Section 6
either to $p=\infty$ or $p \leq 2$.

{\scriptsize \em Acknowledgements: The author would like to thank
Prof. Ping Zhang for the privilege of visiting Beijing and
delivering these lectures, for the kindness and warmth of the
hospitality he has encountered there, and for the generous support
of the Morningside Center. The author would also like to thank his
collaborators Dragos Iftimie, Jim Kelliher, Anna Mazzucato, Helena
Nussenzveig Lopes and Franck Sueur, for the many discussions and
their insight with respect to this material. The author`s research
is supported in part by CNPq grant \# 302.102/2004-3 and by FAPESP
grant \# 07/51490-7.}

\vspace{0.5cm}

\noindent {\sc
Milton C. Lopes Filho\\
Departamento de Matematica, IMECC-UNICAMP.\\
Caixa Postal 6065, Campinas, SP 13083-970, Brasil
\\}
{\it E-mail address:} mlopes@@ime.unicamp.br


\begin{thebibliography}{99}

\bibitem{something} G. K. Batchelor, {\it An introduction to fluid
dynamics} Second edition. Cambridge Mathematical Library. Cambridge
University Press, Cambridge, 1999.

\bibitem{chemin96} J.-Y. Chemin, {\it A remark on the inviscid limit for
two-dimensional incompressible fluids}, Comm. Partial Differential
Equations {\bf 21} (1996), 1771--1779.

\bibitem{CMR98} T. Clopeau, A. Mikelic and R. Robert, {\it On the vanishing viscosity limit
for the 2D incompressible Navier-Stokes equations with friction-type
boundary conditions}, Nonlinearity, {\bf 11} (1998), 1625--1636.

\bibitem{dellelis} C. de Lellis, and L. Szekelyhidi Jr., {\it The Euler equations as a differential
inclusion}, preprint 2007.

\bibitem{E00} Weinan E, {\it Boundary layer theory and the
zero-viscosity limit of the Navier-Stokes equation}, Acta Math.
Sinica, English series, {\bf 16}(2000), 207--218.

\bibitem{EE97} Weinan E and B. Engquist, {\it Blow-up of solutions
to the unsteady Prandtl`s equation}, Comm. Pure Appl. Math. {\bf 50}
(1997), 1287--1293.

\bibitem{grenier00} E. Grenier, {\it On the nonlinear instability of
Euler and Prandtl equations}, Comm. Pure Appl. Math., {\bf LIII}
(2000), 1067-1091.

\bibitem{ILN08} D. Iftimie, M. C. Lopes Filho and H. J. Nussenzveig Lopes, {\it
Incompressible flow around a small obstacle and the vanishing viscosity limit},
in preparation, 2008.

\bibitem{is08} D. Iftimie and F. Sueur, {\it Viscous boundary layers
for the Navier Stokes equations with the Navier slip conditions},
preprint.

\bibitem{jm01} W. J\"{a}ger and A. Mikeli\'{c}, {\it On the roughness-induced effective
boundary conditions for an incompressible viscous flow}, J.
Differential Equations {\bf 170} (2001) 96--122.

\bibitem{Kato84} T. Kato, {\it Remarks on zero viscosity limit for nonstationary Navier-Stokes
flows with boundary} In: S. S. Chern (ed.) {\it Seminar on nonlinear
PDE}, MSRI, 1984.

\bibitem{kelliher06} J. Kelliher, {\it Navier-Stokes equations with
Navier boundary condition for a bounded domain in the plane}, SIAM
J. Math. Anal. {\bf 38} (2006), 210--232.

\bibitem{kelliher07} J. Kelliher, {\it On Kato`s condition for
vanishing viscosity}, Indiana Univ. Math. J. {\bf 56} (2007),
1711-1721.

\bibitem{KLN08} J. Kelliher, M. C. Lopes Filho and H. J. Nussenzveig
Lopes, {\it Vanishing viscosity limit for an expanding domain in
space}, in preparation, 2008.

\bibitem{lieberman} G. Lieberman, {\it  Second order parabolic differential equations.}
World Scientific Publishing Co., Inc., River Edge, NJ, 1996.

\bibitem{lions69} J.-L. Lions, {\it Quelques methodes de resolution des problemes
aux limites non lineaires}, Dunod; Gauthier-Villars, Paris 1969.

\bibitem{LCS03} M. Lombardo, M. Cannone and  M. Sammartino, {\it Well-posedness of
the boundary layer equations}, SIAM J. Math. Anal. {\bf 35} (2003),
no. 4, 987--1004.

\bibitem{LMN06} M. Lopes Filho, A. Mazzucato and H. Nussenzveig
Lopes, {\it Weak solutions, renormalized solutions and enstrophy
defects in 2D turbulence}, Arch. Rat. Mech. Anal., {\bf 179} (2006),
353-387.

\bibitem{LMN08} M. Lopes Filho, A. Mazzucato and H. Nussenzveig
Lopes, {\it Vanishing viscosity limit for incompressible flow inside
a rotating circle}, to appear, Phys. D, 2008.

\bibitem{LMNT07} M. Lopes Filho, A. Mazzucato, H. Nussenzveig
Lopes and M. Taylor, {\it Vanishing Viscosity Limits and Boundary
Layers for Circularly Symmetric 2D Flows}, ArXiV preprint
arXiv:0709.2056v1 [math.AP], 2007, submitted for publication.

\bibitem{LNP05} M. Lopes Filho, H. Nussenzveig Lopes and G. Planas, {\it On the inviscid limit
for two-dimensional incompressible flow with Navier friction
condition}, SIAM J. Appl. Math. {\bf 36}, (2005) no. 4, 1130--1141.

\bibitem{LNX01} M. Lopes Filho, H. Nussenzveig Lopes and
Zhouping Xin, {\it Existence of vortex sheets with reflection
symmetry in two space dimensions} Arch. Rat. Mech. Anal. {\bf 158}
(2001), 235--257.

\bibitem{Majda93} A. Majda, {\it Remarks on weak solutions for vortex sheets with a distinguished
sign}, Indiana Univ. Math. J. {\bf 42} (1993), 921--939.

\bibitem{matsui94} S. Matsui, {\it Example of zero viscosity limit for two-dimensional nonstationary
Navier-Stokes flows with boundary}, Japan J. Indust. Appl. Math.
{\bf 11} (1994),  no. 1, 155--170.

\bibitem{oleinik66} O. Oleinik, {\it On the mathematical theory of boundary layer for an
unsteady flow of incompressible fluid} translated in J. Appl. Math.
Mech. {\bf 30} (1967) 951--974.

\bibitem{oleinik} O. Oleinik, {\it Construction of the solutions of a system of boundary layer
equations by the method of straight lines}, Physics Dokl. 12(1967),
525--528.

\bibitem{oleinikbook} O. Oleinik and V. Samokhin, {\it Mathematical models in
boundary layer theory} Applied Mathematics and Mathematical
Computation, v.15, Chapman \& Hall/CRC, Boca Raton, FL, 1999.

\bibitem{SC98} M. Sammartino and R. Caflisch, {\it Zero viscosity
limit for analytic solutions of the Navier-Stokes equations on a
half-space I and II}, Comm. Math Phys. {\bf 192} (1998), 433-461 and
463-491.

\bibitem{Schlichting} H. Schlichting and K. Gersten, {\it Boundary layer theory}, 8th ed.,
Springer Verlag, Berlin, 2000.

\bibitem{TW97} R. Temam and Xiaoming Wang,
{\it The convergence of the solutions of the Navier-Stokes equations
to that of the Euler equations,} Appl. Math. Lett. {\bf 10} (1997),
29--33.

\bibitem{TW98} R. Temam and Xiaoming Wang, {\it On the behavior of the solutions of the Navier-Stokes
equations at vanishing viscosity}, Annali della Scuola Norm. Sup. Pisa, vol. dedicated to the memory
of E. De Giorgi, Serie IV {\bf XXV} (1998) 807--828.

\bibitem{TW02} R. Temam and Xiaoming Wang, {\it Boundary layer associated with the incompressible
Navier-Stokes equations: the non-characteristic boundary case}, J.
Diff. Eqs. {\bf 179} (2002) 647--686.

\bibitem{Wang01} Xiaoming Wang, {\it A Kato type theorem on zero viscosity limit of Navier-Stokes flows},
Indiana Univ. Math. J. {\bf 50} (2001) no. 1 223--241.

\bibitem{XZ04} Zhouping Xin and Liqun Zhang, {\it On the global existence of solutions to the Prandtl's
system}, Adv. Math. {\bf 181} (2004), no. 1, 88--133.

\bibitem{zz06} Jianwen Zhang and Junning Zhao, {\it On the global
existence and uniqueness of solutions to the nonstationary boundary
layer system}, Science in China, Ser. A Mathematics, {\bf 49}
(2006), 932--960.

\end{thebibliography}
\end{document}